\documentclass[12pt]{amsart}
\usepackage{amssymb}
\usepackage{amsfonts}
\usepackage{latexsym}
\usepackage[all]{xy}
\usepackage{amscd}

\newtheorem{theorem}{Theorem}[section]
\newtheorem{lemma}[theorem]{Lemma}
\newtheorem{proposition}[theorem]{Proposition}
\newtheorem{corollary}[theorem]{Corollary}
\theoremstyle{definition}
\newtheorem{definition}[theorem]{Definition}
\newtheorem{example}[theorem]{Example}

\newtheorem{remark}[theorem]{Remark}
\newtheorem{note}[theorem]{Note}

\newcommand{\C}{\mathcal{C}}
\newcommand{\D}{\mathcal{D}}

\newcommand{\Z}{\mathcal{Z}}
\newcommand{\M}{\mathcal{M}}

\renewcommand{\P}{\mathcal{P}}
\renewcommand{\L}{\mathcal{L}}

\newcommand{\id}{\text{id}}

\newcommand{\ot}{\otimes}
\newcommand{\be}{\mathbf{1}}
\newcommand{\rev}{\text{rev}}

\newcommand{\Hom}{\mbox{Hom}}

\newcommand{\Irr}{\text{Irr}}
\newcommand{\Lagr}{\mbox{Lagr}}

\newcommand{\Rep}{\mbox{Rep}}

\renewcommand{\dim}{\mbox{dim}}
\renewcommand{\deg}{\text{deg }}
\renewcommand{\Vec}{\mbox{Vec}}

\newcommand{\lp}{\left(}
\newcommand{\rp}{\right)}

\newcommand{\la}{\langle}
\newcommand{\ra}{\rangle}

\begin{document}

\title[Lagrangian subcategories of twisted quantum doubles]
{Lagrangian subcategories and braided tensor equivalences
of twisted quantum doubles of finite groups}

\author{Deepak Naidu}
\address{Department of Mathematics and Statistics,
University of New Hampshire, Durham, NH 03824, USA}
\email{dnaidu@unh.edu}

\author{Dmitri Nikshych}
\address{Department of Mathematics and Statistics,
University of New Hampshire, Durham, NH 03824, USA}
\email{nikshych@math.unh.edu}

\date{June 28, 2007}

\begin{abstract}
We classify Lagrangian subcategories of 
the representation category of a twisted quantum double 
$D^\omega(G)$, where $G$ is a finite group and $\omega$ 
is a $3$-cocycle on it. In view of results of \cite{DGNO}
this gives a complete description of all braided tensor
equivalent pairs  of  twisted quantum doubles of finite groups.
We also establish a canonical bijection between 
Lagrangian subcategories of $\Rep(D^\omega(G))$ and
module categories over the category $\Vec_G^\omega$ of twisted 
$G$-graded vector spaces such that the dual tensor category is pointed. 
This can be viewed as  a quantum version of V.~Drinfeld's
characterization of homogeneous spaces of  a Poisson-Lie group
in terms of Lagrangian subalgebras of the double of its   
Lie bialgebra \cite{D}.
As a consequence, we obtain that two group-theoretical fusion categories
are weakly Morita equivalent if and only if their centers
are equivalent as braided tensor categories.
\end{abstract}

\maketitle

\begin{section}
{Introduction}
Throughout this paper we will work over an algebraically closed
field $k$ of characteristic zero. Unless otherwise stated
all cocycles appearing in this 
work will have coefficients in the trivial module $k^\times$.
All categories considered in this work
are assumed to be $k$-linear and semisimple with finite-dimensional
$\Hom$-spaces and finitely many isomorphism classes of simple objects.
All functors are assumed to be additive and $k$-linear.

Let $G$ be a finite group and $\omega$ be a $3$-cocycle on $G$.
In \cite{DPR1, DPR2} R.~Dijkgraaf, 
V.~Pasquier, and P.~Roche introduced a quasi-triangular quasi-Hopf algebra 
$D^\omega(G)$. When $\omega =1$ this quasi-Hopf algebra coincides with 
the Drinfeld double $D(G)$ of $G$ and so $D^\omega(G)$ is often called
a {\em twisted quantum double} of $G$. It is well known that the 
representation category $\Rep(D^\omega(G))$ of $D^\omega(G)$ 
is a modular category  \cite{BK, T} 
and is braided equivalent to the center \cite{K} 
of the tensor category $\Vec_G^\omega$ of finite-dimensional 
$G$-graded vector spaces with the associativity constraint
defined using $\omega$. The category $\Vec_G^\omega$ is a typical
example of a {\em pointed} fusion category, i.e., a finite semisimple
tensor category in which every simple object is invertible.

In \cite{DGNO} a criterion  for a modular category $\C$ to be 
braided tensor equivalent to the center of a category of the form
$\Vec_G^\omega$ for some finite group $G$ and 
$\omega\in Z^3(G,\, k^\times)$ is given.  Namely, such a braided
equivalence exists if and only if $\C$ contains a Lagrangian
subcategory, i.e., a maximal isotropic subcategory of dimension
$\sqrt{\dim(\C)}$.  More precisely, Lagrangian subcategories
of $\C$ parameterize the classes of braided equivalences
between $\C$ and centers of pointed categories, 
see \cite[Section 4]{DGNO}. Note that any Lagrangian
subcategory of $\C$ is equivalent
(as a symmetric tensor category) to the representation category of 
some group by the result of P.~Deligne \cite{De}.

This means that  a description of Lagrangian subcategories 
of  $\Rep(D^\omega(G))$ for
all groups $G$ and $3$-cocycles $\omega$
is  equivalent to a description of all
braided equivalences between representation categories
of twisted group doubles. Such equivalences for elementary
Abelian and extra special groups were studied in \cite{MN}
and \cite{GMN}. A motivation for such study comes from a 
relation between holomorphic orbifolds in the Rational Conformal
Field Theory and twisted group doubles 
observed in \cite{DVVV}, \cite{DPR2}. 

A complete classification  of Lagrangian subcategories of
$\Rep(D^\omega(G))$ is the principal goal of this paper.

\subsection{Main results}

Let $G$ be a finite group and let $\omega \in Z^3(G,\, k^\times)$
be a $3$-cocycle on $G$.

\begin{theorem}
\label{thm 1}
Lagrangian subcategories of the representation category
of the Drinfeld double $D(G)$ are classified by pairs
$(H, B)$, where $H$ is a normal Abelian subgroup of $G$
and $B$ is an alternating $G$-invariant bicharacter on $H$. 
\end{theorem}

The proof is based on the analysis of modular data
(i.e., the $S$- and $T$-matrices) associated to $D(G)$.

Theorem~\ref{thm 1} gives a simple classification of  Lagrangian
subcategories for the untwisted double $D(G)$.  In the twisted
($\omega \neq 1$) case the notion of a $G$-invariant
bicharacter needs to be twisted as well, 
cf.~Definition~\ref{alt omega def}.

\begin{theorem}
\label{thm 2}
Lagrangian subcategories of the representation category
of the twisted double $D^\omega(G)$ are classified by pairs
$(H, B)$, where $H$ is a normal Abelian subgroup of $G$
such that $\omega|_{H\times H\times H}$ is cohomologically
trivial and $B: H\times H \to k^\times$ is a $G$-invariant
alternating $\omega$-bicharacter in the sense of Definition~\ref{alt omega def}.
\end{theorem}

Note that bicharacters in the statement of Theorem~\ref{thm 2} are
in bijection with equivalence classes of $G$-invariant cochains 
$\mu \in C^2(H,\, k^\times)$ such that $\delta^2 \mu =
\omega|_{H\times H\times H}$, see \eqref{Omega 1}.

Let $\L_{(H, B)}$ denote the Lagrangian subcategory of 
$\Rep(D^\omega(G))$ corresponding to a pair $(H, B)$
in Theorems~\ref{thm 1} and \ref{thm 2}. 
Then there is a group $G'$, defined up to an isomorphism,
such that $\L_{(H, B)}$ is equivalent to $\Rep(G')$ 
as a symmetric tensor category, where the braiding of 
$\Rep(G')$  is the trivial one \cite{De}.
The group $G'$ can be described explicitly in terms
of $G$, $H$, and $B$, see Remark~\ref{G'}. Note that
$G \not\cong G'$ in general.

There is a canonical subcategory $\L_{(\{1\},\, 1)} \cong \Rep(G)$
corresponding to the forgetful functor 
$\Rep(D^\omega(G)) \cong \Z(\Vec_G^\omega) \to \Vec_G^\omega$.
We have $\L_{(\{1\},\, 1)} \cap \L_{(H, B)} \cong \Rep(G/H)$.

In \cite{N} the first named author classified indecomposable
$\Vec_G^\omega$-module categories $\M$ with the property 
that the dual fusion category $(\Vec_G^\omega)^*_\M$ 
is pointed.  
Such module categories $\M$ can be thought of as categorical 
analogues  of homogeneous spaces. The above  property
gives rise to an equivalence relation on the set of pairs
$(G,\, \omega)$, where $G$ is a finite group and 
$\omega \in Z^3(G,\, k^\times)$ with 
\begin{equation}
\label{Deepak's equivalence}
(G,\, \omega) \approx (G',\, \omega') \mbox{ if and only if }
(\Vec_G^\omega)^*_\M \cong \Vec_{G'}^{\omega'} \mbox{ for some } \M.
\end{equation}
In other words, $\Vec_G^\omega$ and $\Vec_{G'}^{\omega'}$
are {\em weakly Morita equivalent} in the sense of M.~M\"uger \cite{M}.

\begin{theorem}
\label{thm 3}
There is a canonical bijection between equivalence classes
of indecomposable $\Vec_G^\omega$-module 
categories  $\M$ with respect to which the dual fusion
category $(\Vec_G^\omega)^*_\M$ 
is pointed and Lagrangian subcategories of $\Rep(D^\omega(G))$.
\end{theorem}

Lagrangian subalgebras of the double of a Lie bialgebra
$\mathfrak{g}$ were used by V.~Drinfeld in \cite{D} to describe
Poisson homogeneous spaces  of the Poisson-Lie group $G$
corresponding to $\mathfrak{g}$. So Theorem~\ref{thm 3}
can perhaps be understood as a quantum version of the
correspondence in \cite{D}. Note that quantization of Poisson homogeneous 
spaces was studied in \cite{EK}. 

Recall that a fusion category is called {\em group-theoretical}
if it has a pointed dual.

\begin{theorem}
\label{thm 4}
Let $\C_1$ and $\C_2$ be group-theoretical fusion categories.
Then $\C_1$ and $\C_2$ are weakly Morita equivalent if and
only if their centers $\Z(\C_1)$ and $\Z(\C_2)$ are braided equivalent.
\end{theorem}

If $\C_1$ and $\C_2$ are arbitrary
(i.e., not necessarily group-theoretical) 
weakly Morita equivalent finite tensor categories 
then it was observed by M.~M\"uger in \cite[Remark 3.18]{M} that
$\Z(\C_1)$ is braided tensor equivalent  to $\Z(\C_2)$.
Also, $\Z(\C_1)$ being equivalent  to $\Z(\C_2)$ 
(even in a non-braided way) implies that $\C_1\boxtimes \C_1^\rev$ 
is weakly Morita equivalent to $\C_2\boxtimes \C_2^\rev$,
where $\C^\rev$ denotes the fusion category obtained by reversing
the tensor product in $\C$. At the moment of writing we do not
know if braided equivalence of centers implies weak Morita
equivalence of $\C_1$ and $\C_2$ in a general case.

Note that it was shown by S.~Natale \cite{Nat} that a fusion category
$\C$ is group-theoretical if and only if $\Z(\C)$ is braided  
equivalent to the representation category of some twisted  
group double $D^\omega(G)$ (this  also follows from \cite{O2}).

Combining  the explicit description of weak Morita equivalence classes
of pointed categories from \cite{N} and  correspondence
between braided equivalences of centers and Lagrangian
subcategories from \cite{DGNO} one obtains a complete description
of braided equivalences of twisted quantum doubles.

Recall that two finite groups $G_1$ and $G_2$ were called
{\em categorically Morita equivalent} in \cite{N} if $\Vec_{G_1}$
and $\Vec_{G_2}$ are  weakly Morita equivalent.  Let us write
$G_1 \approx G_2$ for such groups. It follows from  Theorem~\ref{thm 4}
that $G_1 \approx G_2$ if and only if   the corresponding Drinfeld doubles 
$D(G_1)$ and $D(G_2)$ have braided tensor equivalent representation 
categories.  By Theorems~\ref{thm 1} and \ref{thm 3} 
groups categorically Morita equivalent to a given  group $G$
correspond to pairs $(H, B)$, where $H$ is a normal Abelian subgroup
of $G$ and $B$ is an alternating bicharacter on $H$.
Note that such pairs $(H, B)$ for which $B$ is, in addition, a {\em nondegenerate} 
bicharacter were used by P.~Etingof and S.~Gelaki in \cite{EG} to describe groups
{\em isocategorical} to $G$, i.e., such groups $G'$ for
which $\Rep(G') \cong \Rep(G)$ as tensor categories. This non-degeneracy
condition in \cite{EG} is the reason why the categorical Morita equivalence
extends isocategorical equivalence.
Indeed, since $\Rep(D(G))$ is determined by the tensor structure of $\Rep(G)$,
it is clear that isocategorical groups are categorically Morita equivalent. 
On the other hand, the first author constructed in \cite{N} examples
of categorically Morita equivalent but non-isocategorical groups.

For an Abelian group $H$ let $\widehat{H}$ denote the group
of linear characters of $H$.

\begin{corollary}
\label{thm 5}
Let $G, G'$ be finite groups, $\omega \in Z^3(G,\, k^\times)$, 
and $\omega' \in Z^3(G',\, k^\times)$. Then the representation
categories of twisted doubles $D^\omega(G)$ and 
$D^{\omega'}(G')$ are equivalent as braided tensor categories 
if and only if $G$ contains a normal Abelian subgroup $H$
such the following conditions are satisfied:
\begin{enumerate}
\item[(1)]  $\omega|_{H \times H \times H}$ is cohomologically trivial,
\item[(2)] there is a $G$-invariant $($see \eqref{Omega 1}$)$ 
$2$-cochain $\mu \in C^2(H, \, k^\times)$ such that
that $\delta^2 \mu = \omega|_{H \times H \times H}$, and
\item[(3)] there is
an isomorphism $a: G' \xrightarrow{\sim} \widehat{H} \rtimes_{\nu} 
(H \backslash G)$ such that $\varpi \circ (a \times a \times a)$ and 
$\omega'$ are  cohomologically equivalent. 
\end{enumerate}
Here $\nu$ is a certain $2$-cocycle in 
$Z^2(H \backslash G, \, \widehat{H})$ that comes from the $G$-invariance
of $\mu$ and $\varpi$ is a certain $3$-cocycle on $\widehat{H} \rtimes_{\nu} 
(H \backslash G)$ that depends on $\nu$ and on the exact sequence
$1 \to H \to G \to H \backslash G \to 1$ $($see \cite[Theorem 5.8]{N}
for precise definitions$)$.
\end{corollary}

Note that in the special case when $\omega =1$ and $\mu$ is a non-degenerate
$G$-invariant alternating $2$-cocycle on $H$, our construction of the ``dual'' group $G'$
in Corollary~\ref{thm 5}
becomes the construction of a group isocategorical to $G$ from \cite{EG}.
This can be seen by comparing \cite[4.2]{N} and \cite[Formula (2)]{EG}.

\subsection{Organization of the paper}
Section 2 contains necessary preliminary information
about fusion categories, module categories, and modular
categories. We also recall definitions and
results from \cite{DGNO} concerning Lagrangian subcategories
of modular categories.

Section 3 (respectively, Section 4) is devoted to classification 
of Lagrangian categories of  the representation
category of the Drinfeld double (respectively, twisted double)
of a finite group. The reason we prefer to treat untwisted
and twisted cases separately is because our constructions
in the former case do not involve rather technical cohomological
computations present in the latter. We feel that the reader might 
get a better understanding of our results by exploring the untwisted
case first. Of course when $\omega =1$ the results of Section~4 reduce 
to those of Section~3. 

The Sections 3 and 4 contain proofs of
our main results stated above. Theorems~\ref{thm 1} - \ref{thm 4} 
and Corollary~\ref{thm 5} correspond to
Theorems~\ref{untwisted bijection}, \ref{bij 1},
\ref{bijection 1}, \ref{main 1}, and Corollary~\ref{main 2}.

Section 5 contains examples in which we compute  Lagrangian
subcategories of Drinfeld doubles of finite symmetry groups. 
Here we also show that the four non-equivalent non-pointed
fusion categories of dimension $8$ with integral dimensions
of objects are pairwise weakly Morita non-equivalent, and
hence their centers are pairwise non-equivalent as braided 
tensor categories.

\subsection{Acknowledgments}
The present paper would not be possible without \cite{DGNO}
and the second named author is happy to thank his collaborators,
Vladimir Drinfeld, Shlomo Gelaki, and Victor Ostrik for many
useful discussions. 
The authors thank Nicolas Andruskiewitsch, Pavel Etingof, 
and Leonid Vainerman  for helpful comments.
The second author thanks the  
Universit\'e de Caen for hospitality and excellent working conditions
during his visit. The authors were supported by the NSF grant
DMS-0200202. The research of Dmitri Nikshych was supported
by the NSA grant H98230-07-1-0081.

\end{section}
\begin{section}{Preliminaries}

\begin{subsection}{Fusion categories and their module categories}
A {\em fusion category} over $k$ is a
$k$-linear semisimple rigid tensor category with finitely many
isomorphism classes of simple objects, finite-dimensional Hom-spaces,
and simple  neutral object. 

In this paper we only consider fusion categories 
with integral Frobenius-Perron dimensions of simple objects.
It was shown in \cite[Propositions 8.23, 8.24]{ENO} that any such category 
is equivalent to the representation category
of a semisimple quasi-Hopf algebra and 
has a canonical spherical structure with respect
to which the categorical dimension of any object is equal to
its Frobenius-Perron dimension. In particular,
all categorical dimensions are positive integers.
For any object $X$  of a fusion category $\C$ 
let $d(X)$ denote its (Frobenius-Perron) dimension.

By a fusion subcategory of a fusion category we will always mean 
a full fusion subcategory.

A fusion category is said to be {\em pointed} if all its simple
objects are invertible. A typical example of a pointed category 
is $\Vec_G^{\omega}$ - the category of finite-dimensional vector
spaces over $k$ graded by the finite group $G$. The morphisms in this category 
are linear transformations that respect the grading and the associativity
constraint is given by a $3$-cocycle $\omega$ on $G$.

Let $\C = (\C, \, \otimes, \, 1_{\C}, \, \alpha, \, \lambda, \, \rho)$ be a
tensor category, where $1_{\C}$, $\alpha$, $\lambda$, and $\rho$ 
are the unit object, the associativity constraint,
the left unit constraint, and the right unit constraint, respectively.  
A right {\em module category} over $\C$ is 
a category $\M$ together with an exact bifunctor $\otimes: \M \times \C \to 
\M$ and natural isomorphisms
$\mu_{M, \, X, \, Y}: M \otimes (X \otimes Y) \to (M \otimes X) \otimes Y, \,\,
\tau_M: M \otimes 1_{\C} \to M$, for all $M \in \M, \, X, Y \in \C$ such that 
the following two equations hold for all $M \in \M, \, X, Y, Z \in \C$:
\begin{equation*}
\label{module pentagon}
\mu_{M \otimes X, \, Y, \, Z} \, \circ 
\, \mu_{M, \, X, \,  Y \otimes Z} \, \circ 
\, (\id_M \otimes \alpha_{X,Y,Z})
= (\mu_{M, \, X, \, Y}\otimes \id_Z) \, \circ 
\, \mu_{M, \, X \otimes Y, \, Z},
\end{equation*}
\begin{equation*}
\label{module triangle}
(\tau_M \otimes \id_Y) \, \circ
\, \mu_{M, \, 1_\C, \, Y}
= \id_M \otimes \lambda_Y.
\end{equation*}
\par 
Note  that having a $\C$-module structure on a category $\M$
is the same as having a tensor functor from $\C$ to the (strict)
tensor category of endofunctors of $\M$. The coherence conditions on 
a module action follow automatically from those of a tensor functor.

Let $(\M_1, \, \mu^1, \tau^1)$ and $(\M_2, \, \mu^2, \tau^2)$ be two
right module categories over
$\C$. A $\C$-{\em module functor} from $\M_1$ to $\M_2$
is a functor $F: \M_1\to \M_2$ together with natural isomorphisms
$\gamma_{M, \, X}: F(M \otimes X) \to F(M) \otimes X$, for all
$M \in \M_1, \, X \in \C$ such that the following two 
equations hold for all $M \in \M_1, \, X, Y \in \C$:
\begin{equation*}
\label{module functor pentagon}
(\gamma_{M, \, X} \otimes \id_Y) \, \circ 
\, \gamma_{M\otimes X, \, Y} \, \circ
\, F(\mu^1_{M, \, X, \, Y})
= \mu^2_{F(M), \, X, \, Y} \, \circ
\, \gamma_{M, X \otimes Y},
\end{equation*}
\begin{equation*}
\label{module functor triangle}
\tau^1_{F(M)} \, \circ \, \gamma_{M, \, 1_\C} = F(\tau^1_M).
\end{equation*}
Two module categories $\M_1$ and $\M_2$ over $\C$ are {\em equivalent}
if there exists a module functor from $\M_1$ to $\M_2$ which is an
equivalence of categories.
For two module categories $\M_1$ and $\M_2$ over a tensor category
$\C$ their {\em direct sum} is the category $\M_1 \oplus \M_2$ with 
the obvious module category structure. A module category is 
{\em indecomposable} if it is not equivalent to a
direct sum of two non-trivial module categories.

\begin{example}{\bf Indecomposable module categories over pointed categories}.
Let $G$ be a finite group and $\omega \in Z^3(G, \, k^\times)$.
Indecomposable right module categories over $\Vec_G^{\omega}$ correspond to pairs
$(H, \, \mu)$, where $H$ is a subgroup of $G$ such that
$\omega|_{H \times H \times H}$ is cohomologically trivial
and $\mu \in C^2(H, \, k^\times)$ is a $2$-cochain satisfying 
$\delta^2\mu = \omega|_{H \times H \times H}$, i.e.,
\begin{equation}
\label{mu omega}
\mu(h_2, \, h_3)
\mu(h_1h_2, \, h_3)^{-1} \mu(h_1, \, h_2h_3) 
\mu(h_1, \, h_2)^{-1} = \omega(h_1, \, h_2, \,h_3).
\end{equation}
\noindent for all $h_1, h_2, h_3 \in H$ (see \cite{O1}).
Let $\M := \M(H, \, \mu)$ denote the right module category
constructed from the pair $(H, \, \mu)$. The simple objects of
$\M$ are given by the set $H \backslash G$ of right cosets of $H$ in $G$, 
the action of $\Vec_G^{\omega}$
on $\M$ comes from the action of $G$ on $H \backslash G$, and the
module category structure isomorphisms are induced from the $2$-cochain
$\mu$. 
Let $H, H'$ be subgroups of $G$  such that restrictions
of $\omega$ are trivial in $H^3(H, \, k^\times)$ and $H^3(H', \, k^\times)$.
Two pairs $(H, \, \mu)$ and $(H', \, \mu')$, where
$\delta^2 \mu = \omega|_{H \times H \times H}$ and
$\delta^2 \mu' = \omega|_{H' \times H' \times H'}$
give rise to equivalent
$\Vec_G^{\omega}$-module categories if and only if there is $g\in G$
such that $H'=gHg^{-1}$ and $\mu$ and the $g$-conjugate of $\mu'$
differ by a coboundary.
We will say that two elements of $\{ \mu \in C^2(H, \, k^\times) \, | \,
\delta^2 \mu = \omega|_{H \times H \times H}\}$
are equivalent if they differ by a coboundary. Let
\begin{equation}
\label{Omega}
\Omega_{H, \omega} := \text{equivalence classes of }
\left\{ \mu \in C^2(H, \, k^\times)  \mid 
\delta^2 \mu = \omega|_{H \times H \times H}\right\}.
\end{equation}
There is an (in general, 
non-canonical) bijection between $\Omega_{H, \omega}$ and
$H^2(H, \, k^\times)$, i.e., $\Omega_{H, \omega}$ is a 
(non-empty) torsor over $H^2(H, \, k^\times)$. Note that 
$\Omega_{H, 1} = H^2(H, \, k^\times)$.
\end{example}

Let $\M_1$ and $\M_2$ be two right module categories over
a tensor category $\C$. Let $(F^1, \, \gamma^1)$ and $(F^2, \, \gamma^2)$
be module functors from $\M_1$ to $\M_2$. A {\em natural module
transformation} from $(F^1, \, \gamma^1)$ to $(F^2, \, \gamma^2)$ is a 
natural transformation $\eta: F^1 \to F^2$ such that the following
equation holds for all $M \in \M_1$, $X \in \C$:
\begin{equation*}
\label{module trans square}
(\eta_M \otimes \id_X) \, \circ \, \gamma_{M, \, X}^1
= \gamma_{M, \, X}^2 \, \circ \, \eta_{M \otimes X}.
\end{equation*}

Let $\C$ be a tensor category and let $\M$ be a 
right module category over $\C$. 
The {\em dual category} of $\C$ with respect to
$\M$ is  the category $\C^*_\M:=Fun_\C(\M,\M)$ whose objects are $\C$-module 
functors from $\M$ to itself and morphisms are natural 
module transformations.
The category $\C^*_\M$ is a tensor category with tensor product being composition
of module functors.
It is known that if $\C$ is a fusion category
and $\M$ is semisimple $k$-linear and indecomposable,
 then $\C^*_\M$ is a fusion category \cite{ENO}. 

Let $G$ be a finite group and $\omega \in Z^3(G, \, k^\times)$.
For each $x \in G$, 
define $\Upsilon_x : G \times G \to k^\times$ by
\begin{equation}
\label{Upsilon}
\Upsilon_x(g_1, \, g_2) := \frac{\omega(xg_1x^{-1}, \, xg_2x^{-1}, \, x)
\omega(x, \, g_1, \, g_2)}{\omega(xg_1x^{-1}, \, x, \, g_2)}, \qquad 
\mbox{for all } g_1, g_2 \in G.
\end{equation}
\noindent It is straightforward to verify that 
$\delta^2 \Upsilon_x = \frac{\omega}{\omega^x}$, for all
$x \in G$, where 
$$
\omega^x(g_1, \, g_2, \, g_3)
= \omega(xg_1x^{-1}, \, xg_2x^{-1}, \, xg_3x^{-1}),
$$ 
for all $g_1, g_2, g_3 \in G$.

For each $x \in G$, 
define $\nu_x : G \times G \to k^\times$ by
\begin{equation*}
\nu_x(g_1, \, g_2) := \frac{\omega(g_1, \, g_2, \, x) 
\omega(g_1g_2xg_2^{-1}g_1^{-1}, \, g_1, \, g_2)}
{\omega(g_1, \, g_2xg_2^{-1}, \, g_2)}, \qquad \mbox{ for all }
g_1, g_2 \in G.
\end{equation*}
\noindent It is easy to verify that the following relation holds:
\begin{equation}
\label{nu upsilon}
\frac{\Upsilon_{x_1x_2}(g_1, \, g_2)}
{\Upsilon_{x_1}(x_2g_1x_2^{-1}, \, x_2g_2x_2^{-1}) 
\Upsilon_{x_2}(g_1, \, g_2)} = 
\frac{\nu_{g_1}(x_1, \, x_2) \nu_{g_2}(x_1, \, x_2)}
{\nu_{g_1g_2}(x_1, \, x_2)}, \quad \mbox{ for all }
x_1, x_2, g_1, g_2 \in G.
\end{equation}

Let $H$ be a normal subgroup of $G$ such that
$\omega|_{H \times H \times H}$ is cohomologically trivial.
For any $x \in G$ and $\mu \in C^2(H, \, k^\times)$ such that $\delta^2 \mu = 
\omega|_{H \times H \times H}$, 
define 
\[
\mu \triangleleft x := 
\mu^x \times \Upsilon_x|_{H \times H}, 
\]
where
$\mu^x(h_1, \, h_2) = \mu(xh_1x^{-1}, \, xh_2x^{-1})$, for all $h_1, h_2 \in H$. 
It is easy to verify that $\delta^2(\mu \triangleleft x) 
= \omega|_{H \times H \times H}$.
This induces an action of $G$
on $\Omega_{H, \omega}$ (defined in \eqref{Omega}). Indeed, that this
is an action follows from \eqref{nu upsilon}. 
Let $(\Omega_{H, \omega})^G$ denote the set of $G$-invariant elements of
$\Omega_{H, \omega}$, i.e.,
\begin{equation}
\label{Omega 1}
(\Omega_{H, \omega})^G := \left\{\mu
\in \Omega_{H, \omega} \,\, \vline \,\,
\frac{\mu^x}{\mu} \times \Upsilon_x|_{H \times H} \text{ is trivial in } H^2(H, \, k^\times),
\text{ for all } x \in G \right\}.
\end{equation}
%
\begin{example}
\label{pointed mod cats}
{\bf Module categories over $\Vec_G^{\omega}$ with pointed duals}.
Let $G$ be a finite group and $\omega \in Z^3(G, \, k^\times)$. 
It is shown in \cite[Theorem 3.4]{N} that 
the set of equivalence classes of indecomposable module categories over
$\Vec_G^{\omega}$ such that the dual is pointed is in bijection with the set of pairs
$(H, \, \mu)$, where $H$ is a normal Abelian subgroup of $G$ such that
$\omega|_{H \times H \times H}$ is cohomologically trivial and
$\mu \in (\Omega_{H, \omega})^G$ (the description in \cite[Theorem 3.4]{N}
is given in somewhat different but equivalent terms).
\end{example}

Two fusion categories $\C$ and $\D$ are said to be {\em weakly Morita
equivalent} if there exists an indecomposable (semisimple $k$-linear) right
module category $\M$ over $\C$ such that the categories $\C^*_{\M}$ and
$\D$ are equivalent as fusion categories.
It was shown by M.~M\"uger \cite{Mu} that this is indeed an equivalence relation.

A fusion category $\C$ is said to be {\em group theoretical} if it is
weakly Morita equivalent to a pointed category.

\end{subsection}

\begin{subsection}{Modular categories and centralizers}

Let $\C$ be a modular fusion category with braiding $c$, 
twist $\theta$, and S-matrix $S$ (see \cite{BK}). 
Let $\D$ be a full (not necessarily tensor) 
subcategory of $\C$. Its dimension is  defined by 
$\dim(\D) := \sum_{X \in \Irr(\D)} d(X)^2$, where $\Irr(\D)$ 
is the set of isomorphism classes of simple objects in $\D$. 
In \cite{M}, M.~M\"uger introduced the notion of
the {\em centralizer} of $\D$ in $\C$ as 
the fusion subcategory  
\[
\D' := \left\{ X \in \C \mid c(Y, \, X) \circ c(X, \, Y)
= \id_{X \otimes Y}, \mbox{ for all } Y \in \D \right\}. 
\]
It was also shown in \cite{M} that if $\D$
is a fusion subcategory then $\D''=\D$ and
\begin{equation}
\label{dimension}
\dim(\D) \cdot \dim(\D') = \dim(\C).
\end{equation}
Following M.~M\"uger,  we will say that two objects 
$X, Y \in \C$ {\em centralize}  each other if 
$$
c(Y, \, X) \circ c(X, \, Y) = \id_{X \otimes Y}.
$$ 
For simple $X$ and $Y$ this condition is equivalent to
$S(X, \, Y) = d(X) d(Y)$ \cite[Corollary 2.14]{M}.

\begin{remark}
\label{square dim}
If $\D$ is a full subcategory of $\C$ such that all objects 
in $\D$ centralize  each other, i.e., $\D\subseteq \D'$
then $\dim(\D)^2 \leq \dim(\C)$.  
Indeed, we have $\dim(\D) \leq \dim(\D')$
and so it follows from \eqref{dimension}
that $\dim(\D)^2 \leq \dim(\C)$. 
In particular, if $\D$ is a symmetric fusion
subcategory of $\C$, then $\dim(\D)^2 \leq \dim(\C)$.
\end{remark}

\begin{lemma}
\label{lag}
Let $\D$ be a full subcategory of $\C$ $($which is not
{\em a priori} assumed to be closed under the tensor product
or duality$)$ such that $\D \subseteq \D'$.
Then the fusion subcategory $\tilde{\D}
\subseteq \C$  generated by $\D$ is symmetric.
\end{lemma}
\begin{proof}
We may assume that $\D$ is closed under taking duals.
Indeed, it follows from \cite[Proposition 2.12]{ENO}
that $X$ centralizes $Y$ if and only if $X$ centralizes $Y^*$
for any two simple objects $X,Y$ in $\C$.

Let $Z_1, Z_2$ be simple objects in $\tilde{\D}$. 
There exist simple objects $X_1, X_2, Y_1, Y_2$ in $\D$
such that $Z_1$ is contained in $X_1\ot Y_1$ and
$Z_2$ is contained in $X_2\ot Y_2$. 
By  \cite[Lemma 2.4 (i)]{M}, it follows that $Z_1$ centralizes
$X_2 \otimes Y_2$, and hence $Z_1, Z_2$ centralize each other. 
\end{proof}

\begin{corollary}
Let $\D$ be a full subcategory of $\C$ 
such that $\D \subseteq \D'$ and $\dim(\D)^2 = \dim(\C)$.
Then $\D$ is a symmetric fusion subcategory.
\end{corollary}
\end{subsection}

\begin{subsection}{Lagrangian subcategories and braided 
equivalences of twisted group doubles}

Let $\C$ be a modular category. Recall that we chose the canonical
spherical twist for $\C$ with respect
to which the categorical dimension of any object of $\C$ is equal to
its Frobenius-Perron dimension.
This is possible by \cite[Proposition 8.23, 8.24]{ENO}.
Let us recall some definitions and results from \cite{DGNO}.

\begin{definition}
A fusion subcategory $\D\subseteq \C$
is said to be {\em isotropic} if the twist of $\C$ restricts to 
identity on $\D$. An isotropic subcategory $\D\subseteq \C$
is said to be {\em Lagrangian} if $(\dim(\D))^2 = \dim(\C)$.
\end{definition}

An isotropic subcategory $\D$ of $\C$ is necessarily symmetric
and its objects have positive categorical dimensions. 
It follows from \cite{De} that there is a (unique up to an 
isomorphism) group $G$ such that $\D \cong \Rep(G)$
as a symmetric fusion category, where $\Rep(G)$
is considered with its trivial braiding.

Consider the set of all braided tensor equivalences 
$F: \C \xrightarrow{\sim} \Z(\P)$, where $\P$
is a pointed fusion category  
and $\Z(\P)$ denotes its center. There is an equivalence
relation on this set defined as follows. We say that
$F_1: \C \xrightarrow{\sim} \Z(\P_1)$ and 
$F_2: \C \xrightarrow{\sim} \Z(\P_2)$ are equivalent
if there exists a tensor equivalence $\iota: \P_1
\xrightarrow{\sim} \P_2$ such that
$\mathcal{F}_2\circ F_2  = \iota\circ \mathcal{F}_1\circ F_1$,
where $\mathcal{F}_i : \Z(\P_i) \to
\P_i, \, i=1,2$, are the canonical forgetful
functors. Let $\mbox{E}(\C)$ be the collection of
equivalence classes of such equivalences.
Informally, $\mbox{E}(\C)$ is the set of all
``different'' braided equivalences between $\C$
and centers of pointed categories, i.e., 
representation categories of twisted group doubles.

Let $\mbox{Lagr}(\C)$  be the set of all
Lagrangian subcategories of $\C$.

In \cite[Theorem 4.5]{DGNO} it was proved that there is a bijection
\begin{equation}
\label{lagr dgno}
f: \mbox{E}(\C)  \xrightarrow{\sim}  \mbox{Lagr}(\C)
\end{equation}
\noindent defined as follows. Note that each braided tensor
equivalence $F: \C\xrightarrow{\sim} \Z(\P)$ 
gives rise to the Lagrangian subcategory $f(F)$ 
of  $\C$ formed by all objects sent to
multiples of the unit object $\be$ under the forgetful functor
$\Z(\P)\to \P$. This subcategory is clearly
the same for all equivalent choices of $F$.

In particular, the center of a fusion category $\D$ contains
a Lagrangian subcategory if and only if $\D$ is group-theoretical \cite{DGNO}.

\end{subsection}
\begin{subsection}{The Schur multiplier of an Abelian group.}
\label{schur abelian}

Let $H$ be a normal Abelian subgroup of a finite group $G$. 
Let $\Lambda^2H$ denote
the Abelian group of alternating bicharacters on $H$, i.e.,
\begin{equation*}
\Lambda^2H := \left\{B: H \times H \to k^\times \,\, \vline \, 
\begin{tabular}{l}
$B(h_1h_2, \, h) = B(h_1, h)B(h_2, \, h),$\\
$B(h, \, h_1h_2) = B(h, h_1)B(h, \, h_2), \text{ and }$\\
$B(h, \, h) = 1, \mbox{ for all } h, h_1, h_2 \in H$
\end{tabular}
\right\}.
\end{equation*}

Let $Z^2(H, \, k^\times)$ be the group of $2$-cocycles on $H$.
Define a homomorphism 
$alt: Z^2(H, \, k^\times) \to \Lambda^2H : \mu \to alt(\mu)$ by 
\begin{equation*}
alt(\mu)(h_1, \, h_2) := 
\frac{\mu(h_2, \, h_1)}{\mu(h_1, h_2)}, \quad h_1, h_2 \in H.
\end{equation*}

It is well known that $alt$ induces an
isomorphism between the Schur multiplier $H^2(H, \, k^\times)$ of $H$ and  
$\Lambda^2H$. By abuse of notation we denote
this isomorphism also by $alt$:
\begin{equation}
\label{alt}
alt: H^2(H, \, k^\times) 
\xrightarrow{\sim} \Lambda^2H.
\end{equation}

Note that both $H^2(H, \, k^\times)$ and $\Lambda^2H$
are right $G$-modules via the conjugation and that
$alt$ is $G$-linear. 
\end{subsection}
\end{section}
\begin{section}
{Lagrangian subcategories in the untwisted case}
\label{section 3}

We fix notation for this Section.
Let $G$ be a finite group. 
For any $g \in G$, let $K_g$ denote the conjugacy
class of $G$ containing $g$. 
Let $R$ denote a complete set of representatives
of conjugacy classes of $G$. Let $\C$ denote the representation category
$\Rep(D(G))$ of the Drinfeld double of the group $G$:
\begin{equation*}
\C := \Rep(D(G)).
\end{equation*}
The category $\C$ is equivalent to $\Z(\Vec_G)$,
the center of $\Vec_G$. It is well known
that $\C$ is a modular category. 
Let $\Gamma$ denote a complete set of representatives of simple
objects of $\C$. The set $\Gamma$ is in bijection with the set
$\{(a, \, \chi) \mid a \in R \mbox{ and } \chi \mbox{ is an 
irreducible character of }
C_G(a) \}$,
where $C_G(a)$ is the centralizer of $a$ in $G$ (see \cite{CGR}). In what follows we
will identify $\Gamma$ with the previous set,
\begin{equation}
\label{Gamma}
\Gamma := \{(a, \, \chi) \mid a \in R \mbox{ and } \chi \mbox{ is an 
irreducible character of } C_G(a) \}.
\end{equation}
Let $S$ and $\theta$ be (see, e.g. \cite{BK}, \cite{CGR}) 
the $S$-matrix and twist, respectively,
of $\C$. Recall that we take the canonical twist.
It is known that the entries of the $S$-matrix lie in a cyclotomic field.
Also, the values of characters of a finite
group are sums of roots of unity.
So we may assume that all scalars appearing herein are complex numbers;
in particular, complex conjugation and absolute values make sense. 
We have the following formulas for the $S$-matrix, twist and
dimensions:
\begin{equation*}
\begin{split}
S((a, \, \chi), \, (b, \, \chi^\prime)) 
&= \frac{|G|}{|C_G(a)||C_G(b)|} \sum_{g \in G(a, \, b)}
\overline{\chi}(gbg^{-1}) \, \overline{\chi}^\prime(g^{-1}ag),\\
\theta(a, \, \chi) 
& = \frac{\chi(a)}{\deg \chi},\\
d((a, \, \chi))
& = |K_a| \, \deg \chi = \frac{|G|}{|C_G(a)|} \, \deg \chi,
\end{split}
\end{equation*}
\noindent for all $(a, \, \chi), (b, \, \chi^\prime) \in \Gamma$,
where $G(a, \, b) = \{g \in G \mid agbg^{-1} = gbg^{-1}a\}$.

\begin{subsection}{Classification of Lagrangian subcategories
of $\mathbf{\Rep(D(G))}$}

\begin{lemma}
\label{centralize}
Two objects $(a, \, \chi), (b, \, \chi^\prime) \in \Gamma$ centralize each other
if and only if the following conditions hold:\\
(i) The conjugacy classes $K_a, K_b$ commute element-wise,\\
(ii) $\chi(gbg^{-1}) \, \chi^\prime(g^{-1}ag) = \deg \chi \, \deg \chi^\prime$, for all 
$g \in G$.
\end{lemma}
\begin{proof}
By \cite[Corollary 2.14]{M}
two objects $(a, \, \chi), (b, \, \chi^\prime) \in \Gamma$ centralize each other
if and only if 
$$
S((a, \, \chi), (b, \, \chi^\prime)) = \deg \chi \, \deg \chi^\prime.
$$
This is equivalent to the equation
\begin{equation}
\label{eqn}
\sum_{g \in G(a, \, b)}
\chi(gbg^{-1}) \, \chi^\prime(g^{-1}ag) = |G| \, \deg \chi \, \deg \chi^\prime,
\end{equation}
\noindent where $G(a, \, b) = \{g \in G \mid agbg^{-1} = gbg^{-1}a\}$.
It is clear that if the two conditions of the Lemma hold, 
then \eqref{eqn} holds since $G(a, \, b) = G$. 

Now suppose that \eqref{eqn} holds. 
We will show that this implies the two conditions in the statement of
the Lemma. We have
\begin{equation*}
\begin{split}
|G| \, \deg \chi \, \deg \chi^\prime 
&= |\sum_{g \in G(a, \, b)} \chi(gbg^{-1}) \, \chi^\prime(g^{-1}ag)|\\
&\leq \sum_{g \in G(a, \, b)} |\chi(gbg^{-1})| \, |\chi^\prime(g^{-1}ag)|\\
&\leq |G| \, \deg \chi \, \deg \chi^\prime.
\end{split}
\end{equation*}
\noindent So
$\sum_{g \in G(a, \, b)} |\chi(gbg^{-1})| \, |\chi^\prime(g^{-1}ag)|
= |G| \, \deg \chi \, \deg \chi^\prime$. Since 
$$
|G(a, \, b)| \leq |G|,\\
\,\,|\chi(gbg^{-1})| \leq \deg \chi, \mbox{ and } 
|\chi^\prime(g^{-1}ag)| \leq \deg \chi^\prime, 
$$
we must have
$G(a, \, b) = G$, 
$|\chi(gbg^{-1})| = \deg \chi$, and
$|\chi^\prime(g^{-1}ag)| = \deg \chi^\prime$. 
The equality $G(a, \, b) = G$
implies that the conjugacy classes $K_a, K_b$ commute element-wise, 
which is Condition (i) in the
statement of the Lemma. Since
$|\chi(gbg^{-1})| = \deg \chi$, and
$|\chi^\prime(g^{-1}ag)| = \deg \chi^\prime$, there exist roots of
unity $\alpha_g$ and $\beta_g$ such that 
$\chi(gbg^{-1}) = \alpha_g \, \deg \chi$, and
$\chi^\prime(g^{-1}ag) = \beta_g \, \deg \chi^\prime$, for all $g \in G$.
Put this in \eqref{eqn} to get the equation 
\begin{equation}
\label{a}
\sum_{g \in G} \alpha_g \beta_g = |G|.
\end{equation}
\noindent Note that \eqref{a} holds if and only if 
$\alpha_g \beta_g = 1$, for all $g \in G$. 
This is equivalent to
saying that $\chi(gbg^{-1}) \, \chi^\prime(g^{-1}ag) 
= \deg \chi \, \deg \chi^\prime$,
for all $g \in G$ and the Lemma is proved.
\end{proof}

\begin{lemma}
\label{FR}
Let $E$ be a normal subgroup of a finite group $K$. 
Let $\Irr(K)$ denote the set of irreducible characters of $K$.
Let $\rho$ be a $K$-invariant character of $E$ of degree 1. Then
$$
\sum_{\chi \in \Irr(K) : \chi|_E = (\deg \chi) \, \rho} 
(\deg \chi)^2 = \frac{|K|}{|E|}.
$$
\end{lemma}
\begin{proof}
Suppose $\chi$ is any irreducible character of 
$K$. 
Since $\rho$ is $K$-invariant, by Clifford's
Theorem, if $\rho$ is an irreducible constituent of $\chi|_E$, then
\[
\label{chiE}
\chi|_E = (\deg \chi) \, \rho.
\]
By Frobenius reciprocity, the multiplicity of
any irreducible $\chi$ in $\mbox{Hind}_{E}^{K}\rho$ 
is equal to the multiplicity
of $\rho$ in $\chi|_E$.  The latter is equal to $\deg \chi$
if $\chi$ satisfies \eqref{chiE} and $0$ otherwise.
Therefore, 
\[
\sum_{\chi \in \Irr(K) : \chi|_E = (\deg \chi) \, \rho} 
(\deg \chi)^2 = \deg \mbox{Ind}_{E}^{K}\rho = \frac{|K|}{|E|},
\]
as required.
\end{proof}

Let $H$ be a normal Abelian subgroup of $G$ and let $B$ be a
$G$-invariant alternating bicharacter on $H$.
Then $H = \bigcup_{a \in H \cap R}K_a$.
Let 
\begin{equation}
\label{L}
\begin{split}
\L_{(H, \, B)} := 
&\text{ full Abelian subcategory of } \C  \text{ generated by } \\
&\left\{(a, \, \chi) \in \Gamma \,\, \vline \,
\begin{tabular}{l}
$a \in H \cap R \text{ and } \chi \text{ is an irreducible character of } C_G(a)$ \\
$\text{ such that } \chi(h) = B(a, \, h) \,\deg \chi, \text{ for all } h \in H$
\end{tabular}
\right\}.
\end{split}
\end{equation}

\begin{proposition}
\label{Proposition L}
The subcategory $\L_{(H, \, B)} \subseteq \Rep(D(G))$ is Lagrangian.
\end{proposition}
\begin{proof}
We have 
\begin{equation*}
\begin{split}
\chi(gbg^{-1}) \, \chi^\prime(g^{-1}ag) 
& = B(a, \, gbg^{-1}) \, \deg \chi \, B(b, \, g^{-1}ag) 
     \, \deg \chi^\prime\\ 
& = B(a, \, gbg^{-1}) \, B(gbg^{-1}, \,a) 
      \, \deg \chi \,\deg \chi^\prime\\  
& = \deg \chi \,\deg \chi^\prime, 
\end{split}
\end{equation*}
\noindent for all $(a, \, \chi), (b, \, \chi^\prime) \in \L_{(H, \, B)} 
\cap \Gamma, g \in G$.
The second equality above is due to $G$-invariance of $B$ and the third
equality holds since $B$ is alternating. 
By Lemma \ref{centralize}, it follows that
objects in $\L_{(H, \, B)}$ centralize each other.

Also, we have
$\theta_{(a, \, \chi)} = \frac{\chi(a)}{\deg \chi}
= \frac{B(a, \, a)}{\deg \chi} \, \deg \chi = 1$, for all $(a, \, \chi) 
\in \L_{(H, \, B)} \cap \Gamma$. Therefore, $\theta|_{\L_{(H, \, B)}} = \id$.

The dimension of $\L_{(H, \, B)}$ is equal to $|G|$. Indeed, 
\begin{equation*}
\begin{split}
\dim(\L_{(H, \, B)})  
&= \sum_{(a, \, \chi) \in \L_{(H, \, B)} \cap \Gamma} d(a, \, \chi)^2 \\
&= \sum_{(a, \, \chi) \in \L_{(H, \, B)} \cap \Gamma} |K_a|^2 \, (\deg \chi)^2\\
&= \sum_{a \in H \cap R} |K_a|^2 \sum_{\chi :
(a, \, \chi) \in \L_{(H, \, B)} \cap \Gamma} (\deg \chi)^2\\
&= \sum_{a \in H \cap R} |K_a|^2 \frac{|C_G(a)|}{|H|}\\ 
&= \frac{|G|}{|H|} \sum_{a \in H \cap R} |K_a|\\
&= |G|.
\end{split}
\end{equation*}
\noindent The fourth equality above is explained as follows.
Fix $a \in H \cap R$. 
Define $\rho : H \to k^\times$ by 
$\rho(h) := B(a, \, h)$.
Observe that $\rho$ is a $C_G(a)$-invariant
character of $H$ of degree $1$ and then apply Lemma \ref{FR}.

It follows from Lemma \ref{lag} that $\L_{(H, \, B)}$ is a 
Lagrangian subcategory of $\Rep(D(G))$ and the Proposition is proved.
\end{proof}

Now, let $\L$ be a Lagrangian subcategory of $\C$. 
So, in particular, the two 
conditions in Lemma \ref{centralize} 
hold for all simple objects in $\L$.
Define
\begin{equation}
\label{H_L}
H_{\L} := \bigcup_{a \in R : (a, \, \chi) \in \L  
\text{ for some }\chi} K_a.
\end{equation}
Note that $H_{\L}$ is a normal Abelian subgroup of $G$. Indeed, that $H_{\L}$
is a subgroup follows from the fact that $\L$ contains the unit object
and is closed under tensor products. The subgroup $H_{\L}$ is normal in $G$
because it is a union of conjugacy classes of $G$. Finally, that $H_{\L}$ is Abelian
follows by Condition (i) of Lemma \ref{centralize}.

For each $a \in H \cap R$, define 
$\xi_a : H_{\L} \to k^\times$ by 
$$
\xi_a(h) := \frac{\chi(h)}{\deg \chi},
$$
for $h \in H_{\L}$, where $\chi$ is any irreducible character of $C_G(a)$
such that $(a, \, \chi) \in \L \cap \Gamma$. To see that this
definition does not depend on the choice of $\chi$, let 
$(a, \, \chi), \, (a, \, \chi^\prime), \, (b, \, \chi^{\prime\prime})
\in \L \cap \Gamma$ and apply Condition (ii) of Lemma \ref{centralize}
to pairs $(a, \, \chi), \, (b, \, \chi^{\prime\prime})$ and
$(a, \, \chi^\prime), \, (b, \, \chi^{\prime\prime})$  to get 
$$
\frac{\chi(gbg^{-1})}{\deg \chi} = \lp \frac{\chi^{\prime\prime}(g^{-1}ag)}
{\deg \chi^{\prime\prime}} \rp ^{-1} \quad \mbox{ and   } \quad
\frac{\chi^\prime(gbg^{-1})}{\deg \chi^\prime} = \lp \frac {\chi^{\prime\prime}(g^{-1}ag)}
{\deg \chi^{\prime\prime}} \rp ^{-1},
$$ 
for all $g \in G$. This implies that $\frac{\chi|_{H_{\L}}}{\deg \chi}
= \frac{\chi^\prime|_{H_{\L}}}{\deg \chi^\prime}$, for any two pairs
$(a, \, \chi), (a, \, \chi^\prime) \in \L \cap \Gamma$.

For any $a, b \in H_{\L} \cap R$, by Condition (ii) of Lemma \ref{centralize},
$\xi_a$ and $\xi_b$ satisfy the equation:
\begin{equation}
\label{eqn1}
\xi_a(gbg^{-1}) = \xi_b(g^{-1}ag)^{-1}, \quad
\mbox{ for all } g \in G.
\end{equation}
Define a map $B_{\L}: H_{\L} \times H_{\L} \to k^\times$ by
\begin{equation}
\label{B_L}
B_{\L}(h_1, \, h_2) := \xi_a(g^{-1}h_2g),
\end{equation}
\noindent where $h_1 = gag^{-1}, g \in G, a \in H_{\L} \cap R$. 

\begin{proposition}
\label{Proposition B_L}
$B_{\L}$ is a well-defined 
$G$-invariant alternating bicharacter on $H_{\L}$.
\end{proposition}
\begin{proof}
First, let us show that $B_{\L}$ is well-defined. 
Suppose $gag^{-1} = kak^{-1}$,
where $a \in H_{\L} \cap R, g, k \in G$. Then
\begin{equation*}
\begin{split}
B_{\L}(gag^{-1}, \, lbl^{-1})
&= \xi_a((g^{-1}l)b(g^{-1}l)^{-1})\\
&= \xi_b((g^{-1}l)^{-1}a(g^{-1}l))^{-1}\\ 
&= \xi_b(l^{-1}(gag^{-1})l)^{-1}\\
&= \xi_b(l^{-1}(kak^{-1})l)^{-1}\\ 
&= \xi_a((l^{-1}k)^{-1}b(l^{-1}k))\\
&= \xi_a(k^{-1}(lbl^{-1})k)\\
&= B_{\L}(kak^{-1}, \, lbl^{-1}),
\end{split}
\end{equation*}
\noindent for all $b \in H_{\L} \cap R, l \in G$. The second and 
the fifth equalities above are due to \eqref{eqn1}.

Let $h_1 = kak^{-1}, h_2 \in H_{\L}, g \in G$, where $a \in H_{\L} \cap R, k \in G$.
Then
\begin{equation*}
\begin{split}
B_{\L}(gh_1g^{-1}, \, gh_2g^{-1})
&= B_{\L}(gkak^{-1}g^{-1}, \, gh_2g^{-1})\\
&= \xi_a((gk)^{-1}(gh_2g^{-1})(gk))\\
&= \xi_a(k^{-1}h_2k)\\
&= B_{\L}(kak^{-1}, \, h_2)\\
&= B_{\L}(h_1, \, h_2).
\end{split}
\end{equation*}
\noindent So, $B_{\L}$ is $G$-invariant.

Now, 
\begin{equation*}
\begin{split}
B_{\L}(gag^{-1}, \, gag^{-1})
&= B_{\L}(a, \, a)\\
&= \xi_a(a)\\
&= \frac{\chi(a)}{\deg \chi}\\
&= \theta_{(a, \, \chi)}\\
&= 1,
\end{split}
\end{equation*}
\noindent for all $a \in H_{\L} \cap R, g \in G$. The first equality above is
due to the $G$-invariance of $B_{\L}$. So
$B_{\L}(h, \, h) = 1$, for all $h \in H_{\L}$.

Also, $B_{\L}(g_1ag_1^{-1}, \, g_2bg_2^{-1}) B_{\L}(g_2bg_2^{-1}, \, g_1ag_1^{-1})
= \xi_a(g_1^{-1}g_2bg_2^{-1}g_1) \xi_b(g_2^{-1}g_1ag_1^{-1}g_2) = 1$, for all
$g_1, g_2 \in G, a, b \in H \cap R$. 
We used \eqref{eqn1} in the last equality.

To see that $B_{\L}$ is a bicharacter, observe first that $\xi_a$ is a 
homomorphism, for all $a \in H_{\L} \cap R$. We have
\begin{equation*}
\begin{split}
B_{\L}(gag^{-1}, \, h_1) \, B_{\L}(gag^{-1}, \, h_2)
&= \xi_a(g^{-1}h_1g) \, \xi_a(g^{-1}h_2g)\\
&= \xi_a(g^{-1}h_1h_2g)\\
&= B_{\L}(gag^{-1}, \, h_1h_2),
\end{split}
\end{equation*}
\noindent for all $a \in H_{\L} \cap R, g \in G, h_1, h_2 \in H_{\L}$. We
conclude that $B_{\L}$ is a $G$-invariant alternating bicharacter on $H_{\L}$
and the Proposition is proved.
\end{proof}

Recall that $\Lagr(\C)$ denotes the set of Lagrangian subcategories
of a modular category $\C$.

\begin{theorem}
\label{untwisted bijection}
Lagrangian subcategories of the representation category
of the Drinfeld double $D(G)$ are classified by pairs
$(H, B)$, where $H$ is a normal Abelian subgroup of $G$
and $B$ is an alternating $G$-invariant bicharacter on $H$.
\end{theorem}
\begin{proof}
Let $\mathcal{E} := \{(H, \, B) \mid H \mbox{ is a normal Abelian subgroup 
of $G$ and } B \in (\Lambda^2H)^G \}$.
Define a map $\Psi : \mathcal{E} \to \Lagr(\C) : (H, \, B) \mapsto \L_{(H, \, B)}$,
where $\C = \Rep(D(G))$ and $\L_{(H, \, B)}$ is defined in ({\ref{L}}). 
It was shown in Proposition \ref{Proposition L}
that $\L_{(H, \, B)}$ is a Lagrangian subcategory.    

To see that $\Psi$ is 
injective pick any $(H, \, B), (H^\prime, \, B^\prime) \in \mathcal{E}$ and
assume that $\Psi((H, \, B)) = \Psi((H^\prime, \, B^\prime))$. So 
in particular we will have $\L_{(H, \, B)} \cap \Gamma 
= \L_{(H^\prime, \, B^\prime)} \cap \Gamma$.
Note that $H = \cup_{(a, \, \chi) \in \L_{(H, \, B)}\cap \Gamma}K_a$ and 
$H^\prime = \cup_{(a, \, \chi) \in \L_{(H^\prime, \, B^\prime)}\cap \Gamma}K_a$.
Since $\L_{(H, \, B)}\cap \Gamma = \L_{(H^\prime, \, B^\prime)}\cap \Gamma$, it follows that
$H = H^\prime$. Also note that for any $(a, \, \chi) \in \L_{(H, \, B)}\cap \Gamma 
= \L_{(H^\prime, \, B^\prime)}\cap \Gamma$, we have $\chi(h) = B(a, \, h) \, \deg \chi 
= B^\prime(a, \, h) \, \deg \chi$, for all $h \in H = H^\prime$.
Since $B, B^\prime$ are $G$-invariant, it follows that
$B = B^\prime$. So $\Psi$ is injective.

To see that $\Psi$ is surjective pick any $\L \in \Lagr(\C)$. Consider the
pair $(H_{\L}, \, B_{\L})$, where $H_{\L}$ and $B_{\L}$ are defined in
(\ref{H_L}) and (\ref{B_L}), respectively. Proposition \ref{Proposition B_L}
showed that $(H_{\L}, \, B_{\L})$ belongs to the set $\mathcal{E}$. We contend that 
$\Psi((H_{\L}, \, B_{\L})) = \L$. It suffices to show that 
$\L \cap \Gamma \subseteq \L_{(H_{\L}, \, B_{\L})}$. But this hold by definition
of $\L_{(H_{\L}, \, B_{\L})}$ and the observation that 
$\frac{\chi|_{H_{\L}}}{\deg \chi}
= \frac{\chi^\prime|_{H_{\L}}}{\deg \chi^\prime}$, for any two pairs
$(a, \, \chi), (a, \, \chi^\prime) \in \L \cap \Gamma$, $a \in H_{\L} \cap R$. 
So $\Psi$ is surjective and the Theorem is proved.
\end{proof}

\end{subsection}
\begin{subsection}{Bijective correspondence between Lagrangian
subcategories and module categories with pointed duals}

Let $\D$ be a fusion category and let $\M$ be an indecomposable
$\D$-module category. There is a canonical braided tensor
equivalence \cite{Mu, EO} 
\begin{equation}
\label{iota}
\iota_\M: \Z(\D)\xrightarrow{\sim}
\Z(\D^*_\M)
\end{equation}
\noindent defined by identifying both centers with
the category of $\D \boxtimes (\D^*_\M)^{\rev}$-module 
endofunctors of $\M$.

Let $f: E(\C) \xrightarrow{\sim} \Lagr(\C)$ be the bijection
between the set of (equivalence classes of) braided tensor 
equivalences between $\C$ and centers of pointed fusion categories 
and the set of Lagrangian subcategories of $\C$ defined in \cite{DGNO},
see \eqref{lagr dgno}.

\begin{theorem}
\label{bijection}
The assignment $\M \mapsto \iota_\M$ restricts to
a bijection between the set of 
equivalence classes of indecomposable $\Vec_G$-module 
categories  $\M$ with respect to which the dual fusion
category $(\Vec_G)^*_\M$ is pointed and $E( \Rep(D(G)))$. 
\end{theorem}
\begin{proof}
Comparing the result of  \cite{N} (see Example~\ref{pointed mod cats})
and Theorem~\ref{untwisted bijection}
and taking into account that the  isomorphism  
$alt:H^2(H, k^\times) \xrightarrow{\sim} (\Lambda^2H)$
is $G$-linear,
we see that the two sets in question  have the same cardinality.
Thus, to prove the theorem it suffices to check that for 
$\M :=\M(H, \mu)$ one has  $f(\iota_\M) \subseteq \L_{(H, \, alt(\mu))}$,
where  $\L_{(H, \, alt(\mu))}$
is the Lagrangian subcategory defined in \eqref{L}.

By definition, $f(\iota_\M)$ consists of all objects 
$Z$ in $\C= \Z(Vec_G)$ (identified with $\Rep(D(G))$)
such that the $\Vec_G$-module endofunctor
$F_Z: \M \to \M : M \mapsto M \ot Z$ 
is isomorphic 
to a multiple of $\id_\M$. Note that here we abuse
notation and write $Z$ for both object of the center
and its forgetful image.

Let us recall the parameterization of simple objects
of $\Z(\Vec_G)$ in \eqref{Gamma}. Suppose that a simple $Z$ corresponds
to the conjugacy class $K_a$ represented by $a\in R$ and
the character afforded by the irreducible representation $\pi:C_G(a)\to GL(V_\pi)$.
Then as a $G$-graded vector space $Z = \oplus_{x\in K_a}\,V_\pi^x$ 
and the permutation isomorphism $c_{g, Z} : g \ot Z \xrightarrow{\sim} 
Z \ot g$ is induced from $\pi$, where we identify simple objects
of $\Vec_G$ with the elements of the group $G$. 

It is clear that $F_Z$ is isomorphic to a multiple of $\id_\M$ 
as an ordinary functor if and only if $K_a \subseteq H$.
Note that this implies that $H \subseteq C_G(a)$.
Note that for every $\Vec_G$-module functor $F:\M \to \M$
the module functor structure on $F$ is completely determined
by the collection of isomorphisms $F(H1 \ot h) \xrightarrow{\sim}
F(H1) \ot h,\, h\in H$, where $H1$ denotes
the trivial coset in $H \backslash G = \Irr(\M)$.

For $F=F_Z$ the latter isomorphism is given by the composition
\[
(H1 \ot h)\ot Z \xrightarrow{\oplus_{x} \mu(h, x)^{-1} \id_{V_\pi^x}}
H1 \ot (h \ot Z) \xrightarrow{ \id_{H1} 
\ot c_{h, Z}} H1 \ot (Z \ot h)
\xrightarrow{\oplus_{x} \mu(x, h) \id_{V_\pi^x}} (H1 \ot Z)\ot h.
\]
The restriction of $c_{h, Z}$ to $h \ot V_{\pi}^a$ is given
by $\pi(h)$, for all $h \in C_G(a)$.
If the above composition equals identity, then 
$\pi(h) = alt(\mu)(a, \, h) \,\, \id_{V_{\pi}}$, for all $h \in H$.
So $Z \in \L_{(H, \, alt(\mu))}$ and, therefore,
$f(\iota_\M) \subseteq \L_{(H, \, alt(\mu))}$, as required.
\end{proof}

\begin{remark}
\label{G'}
Let us explicitly describe the subcategory $\L_{(H, \, B)}$.
By \cite{De} there is a unique up to an isomorphism group $G'$ 
such that $\L_{(H, \, B)} \cong \Rep(G')$ as a symmetric
category. This group $G'$ is precisely the group of invertible
objects in the dual category $(\Vec_G)^*_{\M(H, \mu)}$, where
$\mu \in Z^2(H,\,k^\times)$ is such that $alt(\mu)=B$. 
It was shown in \cite{N} that $G'$ is an extension
\[
0 \to \widehat{H} \to G' \to H\backslash G \to 0,
\]
with the corresponding second cohomology class being the image
of the cohomology class of $\mu$ under the canonical 
homomorphism $H^2(H,\,k^\times)^G \to  H^2(H\backslash G,\,
\widehat{H})$, see \cite{N} for details.
Note that in general $G \not\cong G'$, see Section~\ref{exs}
for examples.
\end{remark}

\end{subsection}

\end{section}
\begin{section}
{Lagrangian subcategories in the twisted case}

In this Section we extend the constructions of the previous Section
when the associativity is given by a $3$-cocycle 
$\omega\in Z^3(G,\,k^\times)$. 
Note that the results of this Section reduce to the results in 
Section~\ref{section 3} when  $\omega \equiv 1$.

For this Section we follow the notation fixed at the beginning of Section~\ref{section 3}.
Let $\omega$ be a normalized $3$-cocycle on $G$, i.e.,
$\omega$ is a map from $G \times G \times G$ to $k^\times$ satisfying:
\begin{equation}
\label{3-cocycle}
\omega(g_2, \, g_3, \, g_4)\omega(g_1, \, g_2g_3, \, g_4)\omega(g_1, \, g_2, \, g_3)
= \omega(g_1g_2, \, g_3, \, g_4)\omega(g_1, \, g_2, \, g_3g_4),\\
\end{equation}
\begin{equation*}
\omega(g, \, 1_G, \, l) = 1,
\end{equation*}
\noindent for all $g, l, g_1, g_2, g_3, g_4 \in G$.

Let $\C$ denote the representation category
$\Rep(D^{\omega}(G))$ of the twisted quantum double of the group $G$
\cite{DPR1, DPR2}:
\begin{equation*}
\C := \Rep(D^{\omega}(G)).
\end{equation*}
The category $\C$ is equivalent to $\Z(\Vec_G^{\omega})$. It is well known
that $\C$ is a modular category.
Replacing $\omega$ by a cohomologous $3$-cocycle we may assume that 
the values of $\omega$ are roots of unity.

For all $a, g, h \in G$, define
\begin{equation}
\label{beta}
\beta_a(h, g) := \omega(a, \, h, \, g) \omega(h, \, h^{-1}ah, \, g)^{-1}
\omega(h, \, g, \, (hg)^{-1}ahg).
\end{equation}
The $\beta_a$'s satisfy the following equation:
\begin{equation}
\label{beta relation}
\beta_a(x, \, y) \beta_a(xy, \, z) = \beta_a(x, \, yz) \beta_{x^{-1}ax}(y, \, z),
\qquad \mbox{for all } x, y, z \in G.
\end{equation}
Observe that the restriction of each $\beta_a$ to the centralizer $C_G(a)$ of
$a$ in $G$ is a normalized $2$-cocycle.
Let $\Gamma$ denote a complete set of representatives of simple
objects of $\C$. The set $\Gamma$ is in bijection with the set
$\{(a, \, \chi) \mid a \in R \mbox{ and } \chi 
\mbox{ is an irreducible $\beta_a$-character of } C_G(a) \}$.
In what follows we will identify $\Gamma$ with the previous set:
\begin{equation}
\label{Gamma 1}
\Gamma := \{(a, \, \chi) \mid a \in R \mbox{ and } \chi \mbox{ is an 
irreducible $\beta_a$-character of } C_G(a) \}.
\end{equation}
Let $S$ and $\theta$ be the $S$-matrix and twist, respectively,
of $\C$. It is known that the entries of the $S$-matrix lie in a cyclotomic field.
Also, the values of $\alpha$-characters of a finite
group are sums of roots of unity, so they are algebraic numbers, 
where $\alpha$ is any $2$-cocycle whose
values are roots of unity. 
So we may assume that all scalars appearing herein are complex numbers;
in particular, complex conjugation and absolute values make sense. 
We have the following formulas for the $S$-matrix, twist, and dimensions
(see \cite{CGR}):\\\\
$S((a, \, \chi), \, (b, \, \chi^\prime))$ \\
\begin{equation*}
\begin{split}
&= \sum_{g \in K_a, g^\prime \in K_b \cap C_G(g)}
\overline{\lp \frac {\beta_a(x, \, g^\prime) \beta_a(xg^\prime, \, x^{-1})
\beta_b(y, \, g) \beta_b(yg, \, y^{-1})} 
{\beta_a(x, \, x^{-1}) \beta_b(y, \, y^{-1})} \rp}
\overline{\chi}(xg^{\prime}x^{-1}) \, \overline{\chi}^\prime(ygy^{-1}),\\
\theta(a, \, \chi) 
& = \frac{\chi(a)}{\deg \chi},\\
d((a, \, \chi))
& = |K_a| \, \deg \chi = \frac{|G|}{|C_G(a)|} \, \deg \chi,
\end{split}
\end{equation*}
\noindent for all $(a, \, \chi), (b, \, \chi^\prime) \in \Gamma$,
where $g = x^{-1}ax, g^\prime = y^{-1}by$.

\begin{subsection}{Classification of Lagrangian subcategories
of $\mathbf{\Rep(D^{\omega}(G))}$}

\begin{remark}
\label{proj. chars.}
Let $\rho:K \to GL(V)$ be a finite-dimensional projective representation
with $2$-cocycle $\alpha$
on the finite group $K$, i.e., $\rho(xy) = \alpha(x,\,y) \rho(x) \rho(y)$,
for all $x,y\in K$.
Let $\chi$ be the projective character
afforded by $\rho$, i.e., $\chi(x) = \text{Trace}(\rho(x))$, for all $x \in K$.
Suppose that the values of $\alpha$ are roots of
unity. Then $|\chi(x)| \leq \deg \chi$, for all
$x \in K$ and we have equality if and only if $\rho(x) \in k^\times \cdot \id_V$.\\
\end{remark} 

\begin{lemma}
\label{centralize 1}
Two objects $(a, \, \chi), (b, \, \chi^\prime) \in \Gamma$ centralize each other
if and only if the following conditions hold:\\
(i) The conjugacy classes $K_a, K_b$ commute element-wise,\\
(ii) $\lp \frac {\beta_a(x, \, y^{-1}by) \beta_a(xy^{-1}by, \, x^{-1})
\beta_b(y, \, x^{-1}ax) \beta_b(yx^{-1}ax, \, y^{-1})}
{\beta_a(x, \, x^{-1}) \beta_b(y, \, y^{-1})} \rp
\chi(xy^{-1}byx^{-1}) \, \chi^\prime(yx^{-1}axy^{-1})
= \deg \chi \, \deg \chi^\prime$, for all 
$x, y \in G$.
\end{lemma}
\begin{proof}
Two objects $(a, \, \chi), (b, \, \chi^\prime) \in \Gamma$ centralize each other
if and only if $S((a, \, \chi), (b, \, \chi^\prime)) = \deg \chi \, \deg \chi^\prime$.
This is equivalent to the equation:
\begin{equation}
\begin{split}
\label{eqn 1}
\sum_{g \in K_a, g^\prime \in K_b \cap C_G(g)}
\lp \frac {\beta_a(x, \, g^\prime) \beta_a(xg^\prime, \, x^{-1})
\beta_b(y, \, g) \beta_b(yg, \, y^{-1})} 
{\beta_a(x, \, x^{-1}) \beta_b(y, \, y^{-1})} \rp
\chi(xg^{\prime}x^{-1}) \, \chi^\prime(ygy^{-1}) \\
 = |K_a||K_b| \, \deg \chi \, \deg \chi^\prime,
\end{split}
\end{equation}
\noindent where $g = x^{-1}ax, g^\prime = y^{-1}by$.
It is clear that if the two conditions of the Lemma hold, then \eqref{eqn 1}
holds since the set over which the above sum is taken is equal to $K_a \times K_b$. 

Now suppose that \eqref{eqn 1}
holds. We will show that this implies the two conditions in the statement of
the Lemma. We have 
\begin{equation*}
\begin{split}
&|K_a||K_b| \, \deg \chi \, \deg \chi^\prime\\
&= \left|\sum_{g \in K_a, g^\prime \in K_b \cap C_G(g)}
\lp \frac {\beta_a(x, \, g^\prime) \beta_a(xg^\prime, \, x^{-1})
\beta_b(y, \, g) \beta_b(yg, \, y^{-1})} 
{\beta_a(x, \, x^{-1}) \beta_b(y, \, y^{-1})} \rp
\chi(xg^{\prime}x^{-1}) \, \chi^\prime(ygy^{-1}) \right|\\
&\leq 
\sum_{g \in K_a, g^\prime \in K_b \cap C_G(g)}
\left|\lp \frac {\beta_a(x, \, g^\prime) \beta_a(xg^\prime, \, x^{-1})
\beta_b(y, \, g) \beta_b(yg, \, y^{-1})} 
{\beta_a(x, \, x^{-1}) \beta_b(y, \, y^{-1})} \rp \right|
\left|\chi(xg^{\prime}x^{-1})\right| \, \left|\chi^\prime(ygy^{-1})\right|\\
&= 
\sum_{g \in K_a, g^\prime \in K_b \cap C_G(g)}
\left|\chi(xg^{\prime}x^{-1})\right| \, \left|\chi^\prime(ygy^{-1})\right|\\
&\leq |K_a||K_b| \, \deg \chi \, \deg \chi^\prime
\end{split}
\end{equation*}

So 
\begin{equation*}
\sum_{g \in K_a, g^\prime \in K_b \cap C_G(g)}
\left|\chi(xg^{\prime}x^{-1})\right| \, |\chi^\prime(ygy^{-1})|
= |K_a||K_b| \, \deg \chi \, \deg \chi^\prime.
\end{equation*}

Since $|\{(g, \, g^\prime) \mid g \in K_a, g^\prime \in K_b \cap C_G(g)\}| 
\leq |K_a||K_b|$,
$|\chi(xg^{\prime}x^{-1})| \leq \deg \chi$, and \\
$|\chi^\prime(ygy^{-1})| \leq \deg \chi^\prime$, we must have
$|\{(g, \, g^\prime) \mid g \in K_a, g^\prime \in K_b \cap C_G(g)\}| = |K_a||K_b|$,
i.e. $\{(g, \, g^\prime) \mid g \in K_a, g^\prime \in K_b \cap C_G(g)\}
= K_a \times K_b$, 
$|\chi(xg^{\prime}x^{-1})| = \deg \chi$, and
$|\chi^\prime(ygy^{-1})| = \deg \chi^\prime$. 
The equality $\{(g, \, g^\prime) \mid g \in K_a, g^\prime \in K_b \cap C_G(g)\}
= K_a \times K_b$
implies that $K_b \subseteq C_G(g)$, for all $g \in K_a$ . This is equivalent to
the condition that $K_a, K_b$ commute element-wise which is Condition (i) in
the statement of the Lemma. 
Now, \eqref{eqn 1} becomes:
\begin{equation}
\label{eqn 2}
\sum_{(g, \, g^\prime) \in K_a \times K_b}
\lp \frac {\beta_a(x, \, g^\prime) \beta_a(xg^\prime, \, x^{-1})
\beta_b(y, \, g) \beta_b(yg, \, y^{-1})} 
{\beta_a(x, \, x^{-1}) \beta_b(y, \, y^{-1})} \rp
\frac {\chi(xg^{\prime}x^{-1})}{\deg \chi} \, 
\frac{\chi^\prime(ygy^{-1})}{\deg \chi^\prime} 
 = |K_a||K_b|, 
\end{equation}
\noindent where $g = x^{-1}ax, g^\prime = y^{-1}by$.
Since $|\chi(xg^{\prime}x^{-1})| = \deg \chi$, and
$|\chi^\prime(ygy^{-1})| = \deg \chi^\prime$, by Remark \ref{proj. chars.},
$\frac {\chi(xg^{\prime}x^{-1})}{\deg \chi}$ and
$\frac{\chi^\prime(ygy^{-1})}{\deg \chi^\prime}$ are roots of unity.
Note that \eqref{eqn 2} holds if and only if
\begin{equation*}
\lp \frac {\beta_a(x, \, g^\prime) \beta_a(xg^\prime, \, x^{-1})
\beta_b(y, \, g) \beta_b(yg, \, y^{-1})} 
{\beta_a(x, \, x^{-1}) \beta_b(y, \, y^{-1})} \rp
\chi(xg^{\prime}x^{-1}) \, \chi^\prime(ygy^{-1}) = \deg \chi \,\,
\deg \chi^\prime,
\end{equation*}
\noindent for all $g \in K_a, g^\prime \in K_b$, where
$g = x^{-1}ax, g^\prime = y^{-1}by$. This is equivalent to 
Condition (ii) in the statement of the Lemma.
\end{proof}

\begin{note}
Let $E$ be a subgroup of a finite group $K$. Let $\alpha$ be a $2$-cocycle on $K$.
Let $\chi$ be a projective $\alpha$-character of $E$.
For any $x \in K$, define $\chi^x$ by 
$$
\chi^x(l) := \alpha(lx, \, x^{-1})^{-1} \alpha(x, \, x^{-1}lx)^{-1}
\alpha(x, \, x^{-1}) \, \chi(x^{-1}lx), 
$$
for all $l \in E$. Then $\chi^x$ is  
a projective $\alpha$-character of $xEx^{-1}$. 
Suppose $E$ is normal in $K$. Then $\chi$ is said to be $K$-{\em invariant} if 
$\chi^x = \chi$, for all $x \in K$.
\end{note}

\begin{lemma}
\label{FR 1}
Let $E$ be a normal subgroup of a finite group $K$. 
Let $\alpha$ be a $2$-cocycle on $K$.
Let $\Irr(K)$ denote the set of irreducible projective
$\alpha$-characters of $K$.
Let $\rho$ be a $K$-invariant projective
$\alpha|_{E\times E}$-character of $E$ of degree $1$. Then
$$
\sum_{\chi \in \Irr(K) : \chi|_E = (\deg \chi) \, \rho} 
(\deg \chi)^2 = \frac{|K|}{|E|}.
$$
\end{lemma}
\begin{proof}
The proof is completely similar to the one given in
Lemma \ref{FR} except in this case we apply Clifford's Theorem
\cite[Theorem 8.1]{KAR}
and Frobenius reciprocity \cite[Proposition 4.8]{KAR} for projective characters.
\end{proof}

Let $H$ be a normal Abelian subgroup of $G$.  

Recall that $\omega \in Z^3(G,\, k^\times)$ gives rise to a collection
\eqref{beta relation}  of $2$-cochains $\beta_a,\, a\in G$.

\begin{definition}
We will say that a map $B: H \times H \to k^\times$ 
is an {\em alternating $\omega$-bicharacter on $H$} if it satisfies the following
three conditions:
\begin{eqnarray}
\label{1 1}
B(h_1, \, h_2) = B(h_2, \, h_1)^{-1},\\
\label{1 2}
B(h, \, h) = 1,\\
\label{1 3}
\delta^1B_h = \beta_{h}|_{H \times H},
\end{eqnarray}
for all $h, h_1, h_2 \in H$, where the map $B_h: H \to k^\times$ is defined by
$B_h(h_1) := B(h, \, h_1)$, for all $h, h_1 \in H$.
\end{definition}

\begin{definition}
\label{alt omega def}
We will say that an alternating $\omega$-bicharacter 
$B: H \times H \to k^\times$ on $H$ is {\em $G$-invariant} if it satisfies the following
condition:
\begin{equation}
\label{1 4}
B(x^{-1}ax, \, h) = \frac{\beta_a(x,\, h) \beta_a(xh, \, x^{-1})}
{\beta_a(x, \, x^{-1})} 
\, B(a, \, xhx^{-1}), \quad \mbox{ for all } 
x \in G, a \in H \cap R, h \in H.
\end{equation}
\end{definition}

Define
\begin{equation}
\begin{split}
\label{alt omega}
\Lambda_{\omega}^2H := \{B: H \times H \to k^\times \mid
B \text{ is an alternating } \omega-\text{bicharacter on } H\},
\end{split}
\end{equation}
and
\begin{equation}
\begin{split}
\label{alt omega 1}
(\Lambda_{\omega}^2H)^G := \{B \in \Lambda_{\omega}^2H \mid 
B \text{ is $G$-invariant} \}.
\end{split}
\end{equation} 

\begin{remark}
If $\omega \equiv 1$, then $(\Lambda_{\omega}^2H)^G$ is the Abelian
group of $G$-invariant alternating bicharacters on $H$.
\end{remark}

\begin{remark}
If $B$ is an alternating $\omega$-bicharacter on $H$, then the
restriction $\omega|_{H \times H \times H}$ must be
cohomologically trivial. Indeed, 
let $\omega_H := \omega|_{H \times H \times H}$. Then 
$B$ defines a braiding on the
fusion category $\Vec_H^{\omega_H}$.  The isomorphism 
$h_1 \ot h_2 \xrightarrow{\sim} h_2 \ot h_1$ is given by $B(h_1, \, h_2)$,
for all $h_1, h_2 \in H$,
where we identify simple objects of $\Vec_H^{\omega_H}$ with elements of $H$.
It is known ( see, e.g., \cite{Q}, \cite{FRS}) that in this case $\omega_H$ 
is an {\em Abelian} $3$-cocycle on $H$. By a classical result
of Eilenberg and MacLane \cite{EM} the third Abelian cohomology group of $H$
is isomorphic to the (multiplicative) group of quadratic forms on $H$.
The value of the corresponding quadratic form $q$ on $h\in H$
is given by $q(h) = B(h,h)$. Since $B$ is alternating we have
$q\equiv 1$ and so $\omega_H$ must be cohomologically trivial.
\end{remark}

Let $B \in (\Lambda^2_{\omega}H)^G$ and define:
\begin{equation}
\label{L 1}
\begin{split}
\L_{(H, \, B)} := 
&\text{ full Abelian subcategory of } \C \text{ generated by } \\
&\left\{(a, \, \chi) \in \Gamma \,\, \vline \,
\begin{tabular}{l}
$a \in H \cap R \text{ and } \chi \text{ is an irreducible $\beta_a$-character of } C_G(a)$ \\
$\text{ such that } \chi(h) = B(a, \, h) \,\deg \chi, \text{ for all } h \in H$
\end{tabular}
\right\}
\end{split}
\end{equation}

\begin{proposition}
\label{Proposition L 1}
The subcategory $\L_{(H, \, B)} \subseteq \Rep(D^{\omega}(G))$ is Lagrangian.
\end{proposition}
\begin{proof}
Pick any $(a, \, \chi), (b, \, \chi^\prime) \in \L_{(H, \, B)}\cap \Gamma$. We have 
\begin{equation*}
\begin{split}
&\lp \frac {\beta_a(x, \, y^{-1}by) \beta_a(xy^{-1}by, \, x^{-1})
\beta_b(y, \, x^{-1}ax) \beta_b(yx^{-1}ax, \, y^{-1})}
{\beta_a(x, \, x^{-1}) \beta_b(y, \, y^{-1})} \rp
\chi(xy^{-1}byx^{-1}) \, \chi^\prime(yx^{-1}axy^{-1})\\
& = \frac {\beta_a(x, \, y^{-1}by) \beta_a(xy^{-1}by, \, x^{-1})}
{\beta_a(x, \, x^{-1})} B(a, \, xy^{-1}byx^{-1}) \\
& \hspace{1.2in} \times \frac {\beta_b(y, \, x^{-1}ax) \beta_b(yx^{-1}ax, \, y^{-1})}
{\beta_b(y, \, y^{-1})} B(b, \, yx^{-1}axy^{-1})
\times \deg \chi  \, \deg \chi^\prime\\ 
& = B(x^{-1}ax, \, y^{-1}by) \, B(y^{-1}by, \, x^{-1}ax) \, \deg \chi
\, \deg \chi^{\prime}\\  
& = \deg \chi \,\deg \chi^\prime, 
\end{split}
\end{equation*}
\noindent for all $x, y \in G$.
The second equality above is due to \eqref{1 4}
while the third equality is due to \eqref{1 1}.
Note that $K_a, K_b$ commute element-wise since $H$ is Abelian.
By Lemma \ref{centralize 1}, it follows that
objects in $\L_{(H, \, B)}$ centralize each other.

Also, $\theta|_{\L_{(H, \, B)}} = \id$. The proof of this assertion
is exactly the one given in Proposition \ref{Proposition L}. 

Now, fix $a \in H \cap R$ and observe that $B_a$ defines a $C_G(a)$-invariant
$\beta_a$-character of $H$ of degree $1$. Indeed, 
\begin{equation*}
\begin{split}
(B_a)^x(h)
&= \frac{\beta_a(x, \, x^{-1})}{\beta_a(hx, \, x^{-1}) \beta_a(x, \, x^{-1}hx)}
B(a, \, x^{-1}hx)\\
&= B(x^{-1}ax, \, x^{-1}hx)^{-1} \, B(a, \, h) \, B(a, \, x^{-1}hx)\\
&= B(a, \, h),
\end{split}
\end{equation*}
for all $x \in C_G(a), h \in H$. The second equality above
is due to \eqref{1 4}.

The dimension of $\L_{(H, \, B)}$ is equal to $|G|$. 
The proof of this assertion is exactly
the one given in Proposition \ref{Proposition L} except we
appeal to Lemma \ref{FR 1} in this case.

It follows from Lemma \ref{lag} that $\L_{(H, \, B)}$ is a 
Lagrangian subcategory of $\Rep(D^{\omega}(G))$ and the Proposition is proved.
\end{proof}

\begin{lemma}
\label{coboundary}
Let $H$ be a normal Abelian subgroup of $G$. 
Let $B:H \times H \to k^\times$ be a map satisfying \eqref{1 1},\eqref{1 2}, and
\eqref{1 4}.
Suppose $\delta^1B_a = \beta_{a}|_{H \times H},
\mbox{ for all } a \in H \cap R$.
Then $B \in (\Lambda^2_{\omega}H)^G$.
\end{lemma}
\begin{proof}
We only need to verify that \eqref{1 3} holds.
We have
\begin{equation*}
\begin{split}
& (\delta^1B_{x^{-1}ax})(h_1, \, h_2) \\
& = \frac{B(x^{-1}ax, \, h_1) B(x^{-1}ax, \, h_2)}{B(x^{-1}ax, \, h_1h_2)}\\
& = \lp \frac{\beta_a(x, \, h_1) \beta_a(xh_1, \, x^{-1})}
{\beta_a(x, \, x^{-1})} \rp B(a, \, xh_1x^{-1}) \times
\lp \frac{\beta_a(x, \, h_2) \beta_a(xh_2, \, x^{-1})}
{\beta_a(x, \, x^{-1})} \rp B(a, \, xh_2x^{-1}) \\
& \hspace{2.3in} \times \lp \frac{\beta_a(x, \, x^{-1})}
{\beta_a(x, \, h_1h_2) \beta_a(xh_1h_2, \, x^{-1})}\rp B(a, \, xh_1h_2x^{-1})^{-1}\\
& = \frac{\beta_a(x, \, h_1) \beta_a(xh_1, \, x^{-1}) 
\beta_a(x, \, h_2) \beta_a(xh_2, \, x^{-1}) \beta_a(xh_1x^{-1}, \, xh_2x^{-1})}
{\beta_a(x, \, x^{-1}) \beta_a(x, \, h_1h_2) \beta_a(xh_1h_2, \, x^{-1})} \\
& = \frac{\beta_{x^{-1}ax}(h_1, \, h_2) \beta_a(xh_1, \, x^{-1}) 
\beta_a(x, \, h_2) \beta_a(xh_2, \, x^{-1}) \beta_a(xh_1x^{-1}, \, xh_2x^{-1})}
{\beta_a(xh_1, \, h_2) \beta_a(x, \, x^{-1}) \beta_a(xh_1h_2, \, x^{-1})}\\
& = \frac{\beta_{x^{-1}ax}(h_1, \, h_2) \beta_a(xh_1, \, h_2x^{-1}) 
\beta_{x^{-1}ax}(x^{-1}, \, xh_2x^{-1}) \beta_a(x, \, h_2x^{-1}) 
\beta_{x^{-1}ax}(h_2, \, x^{-1})}
{\beta_a(xh_1, \, h_2) \beta_a(x, \, x^{-1}) \beta_a(xh_1h_2, \, x^{-1})}\\ 
& = \beta_{x^{-1}ax}(h_1, \, h_2),  
\end{split}
\end{equation*}
\noindent for all $x \in G, a \in H \cap R, h_1, h_2 \in H$. In the
second equality above, we used \eqref{1 4}.
In the third equality we used $\delta^1B_a = \beta_{a}|_{H \times H}$.
and canceled some factors.
In the fourth equality we used \eqref{beta relation} with
$(x, \, y, \, z) = (x, \, h_1, \, h_2)$. In the fifth equality 
we used \eqref{beta relation} twice with
$(x, \, y, \, z) = (x, \, h_2, \, x^{-1}), (xh_1, \, x^{-1}, \, xh_2x^{-1})$.
In the last equality we used \eqref{beta relation} twice with
$(x, \, y, \, z) = (xh_1, \, h_2, \, x^{-1}), (x, \, x^{-1}, \, xh_2x^{-1})$.
\end{proof}

Now, let $\L$ be a Lagrangian subcategory of $\C$. So, in particular, the two 
conditions in Lemma \ref{centralize 1} hold for all objects in $\L \cap \Gamma$.
Define
\begin{equation}
\label{H_L 1}
H_{\L} := \bigcup_{a \in R : (a, \, \chi) \in \L \cap \Gamma 
\text{ for some }\chi} K_a.
\end{equation}
\noindent Note that $H_{\L}$ is a normal Abelian 
subgroup of $G$.

Define a map $B_{\L}: H_{\L} \times H_{\L} \to k^\times$ by
\begin{equation}
\label{B_L 1}
B_{\L}(h_1, \, h_2) := \frac {\beta_a(x, \, h_2) \beta_a(xh_2, \, x^{-1})}
{\beta_a(x, \, x^{-1})} \times \frac{\chi(xh_2x^{-1})}{\deg \chi},
\end{equation}
\noindent where $h_1 = x^{-1}ax, x \in G, a \in H_{\L} \cap R$ and
$\chi$ is any $\beta_a$-character of $C_G(a)$ such that 
$(a, \, \chi) \in \L \cap \Gamma$. The above definition does not
depend on the choice of $\chi$. The proof of this assertion is similar 
to the proof given for the corresponding assertion in the untwisted case.

\begin{proposition}
\label{Proposition B_L 1}
The map $B_{\L}$ defined in \eqref{B_L 1} is an element of
$(\Lambda^2_{\omega}H)^G$.
\end{proposition}
\begin{proof}
First, let us show that $B_{\L}$ is well-defined. 
Suppose $x^{-1}ax = z^{-1}az$,
where $a \in H_{\L} \cap R, x, z \in G$. Then
\begin{equation*}
\begin{split}
B_{\L}(x^{-1}ax, \, y^{-1}by)
&= \frac {\beta_a(x, \, y^{-1}by) \beta_a(xy^{-1}by, \, x^{-1})}
{\beta_a(x, \, x^{-1})} \times \frac{\chi(xy^{-1}byx^{-1})}{\deg \chi}\\
&= \lp \frac {\beta_b(y, \, x^{-1}ax) \beta_b(yx^{-1}ax, \, y^{-1})}
{\beta_b(y, \, y^{-1})} \rp ^{-1}
\lp \frac{\chi^\prime(yx^{-1}axy^{-1})}{\deg \chi^{\prime}} \rp ^{-1}\\ 
&= \lp \frac {\beta_b(y, \, z^{-1}az) \beta_b(yz^{-1}az, \, y^{-1})}
{\beta_b(y, \, y^{-1})} \rp ^{-1}
\lp \frac{\chi^\prime(yz^{-1}azy^{-1})}{\deg \chi^{\prime}} \rp ^{-1}\\
&= \frac {\beta_a(z, \, y^{-1}by) \beta_a(zy^{-1}by, \, z^{-1})}
{\beta_a(z, \, z^{-1})} \times \frac{\chi(zy^{-1}byz^{-1})}{\deg \chi}\\ 
&= B_{\L}(z^{-1}az, \, y^{-1}by),
\end{split}
\end{equation*}
\noindent for all $b \in H_{\L} \cap R, y \in G$,
where $\chi^\prime$ is any irreducible $\beta_b$-character
of $C_G(b)$ such that $(b, \, \chi^\prime) \in \L \cap \Gamma$. 
The second and the fourth equalities above
are due to Condition (ii) of Lemma \ref{centralize 1}.

The map $B_{\L}$ satisfies \eqref{1 1}
because Condition (ii) of Lemma \ref{centralize 1} holds. 
Let us show that \eqref{1 2} holds for $B_{\L}$:
\begin{equation*}
\begin{split}
&B_{\L}(x^{-1}ax, \, x^{-1}ax) \\
& = \frac{\beta_a(x, \, x^{-1}ax) \beta_a(ax, \, x^{-1})}{\beta_a(x, \, x^{-1})}
\times \frac{\chi(a)}{\deg \chi}\\
& = \omega(a, \, x, \, x^{-1}ax) \times
\frac{\omega(a, \, ax, x^{-1}) \omega(ax, \, x^{-1}, \, a)}
{\omega(ax, \, x^{-1}ax, x^{-1})} \times
\frac{\omega(x, \, x^{-1}ax, \, x^{-1})}
{\omega(a, \, x, \, x^{-1}) \omega(x, \, x^{-1}, \, a)} 
\times \theta_{(a, \, \chi)}\\
& = \frac{\omega(a, \, x, \, x^{-1}a) \omega(ax, \, x^{-1}, \, a)}
{\omega(a, \, x, \, x^{-1}) \omega(x, \, x^{-1}, \, a)}\\
& =  1,
\end{split}
\end{equation*}
\noindent for all $x \in G, a \in H_{\L} \cap R$. 
In the second equality we used the definition
of $\beta_a$. In the third equality we used \eqref{3-cocycle}
with $(g_1, \, g_2, \, g_3, \, g_4) = (a, \, x, \, x^{-1}ax,\, x^{-1})$ and
used the fact that $\theta_{(a, \, \chi)} = 1$ . In the fourth equality
we used  \eqref{3-cocycle} with $(g_1, \, g_2, \, g_3, \, g_4) =
(a, \, x, \, x^{-1}, \, a)$.

The map $B_{\L}$ satisfies \eqref{1 4}
because $B_{\L}(a, \, xhx^{-1}) = \frac{\chi(xhx^{-1})}{\deg \chi}$,
for all $a \in H \cap R, x \in G, h \in H$. 
We have $B_{\L}(a, \, h_1) B_{\L}(a, \, h_2) = \frac{\chi(h_1)}{\deg \chi}
\frac{\chi(h_2)}{\deg \chi} = \beta_a(h_1, \, h_2)\frac{\chi(h_1h_2)}{\deg \chi} 
= \beta_a(h_1, \, h_2) \, B_{\L}(a, \, h_1h_2)$, for all $a \in H \cap R, h_1, h_2 \in H$.
The second last equality above is because $H$ acts as scalars on the projective
$\beta_a$-representation of $C_G(a)$ whose projective character is $\chi$.
By Lemma \ref{coboundary} it follows that $B_{\L} \in \Lambda^2_{\omega}H$
and the Proposition is proved.
\end{proof}

\begin{theorem}
\label{bij 1}
Lagrangian subcategories of the representation category
of the twisted double $D^\omega(G)$ are classified by pairs
$(H, B)$, where $H$ is a normal Abelian subgroup of $G$
such that $\omega|_{H\times H\times H}$ is cohomologically
trivial and $B: H\times H \to k^\times$ is a $G$-invariant
alternating $\omega$-bicharacter in the sense of Definition~\ref{alt omega def}.
\end{theorem}
\begin{proof}
The proof is completely similar to the one given in 
Theorem \ref{untwisted bijection}.
\end{proof}

\end{subsection}
\begin{subsection}{Bijective correspondence between Lagrangian
subcategories and module categories with pointed duals}

Recall \cite{N} that equivalence classes of indecomposable 
module categories over $\Vec_G^{\omega}$
for which the dual is pointed are in bijection with pairs $(H, \, \mu)$,
where $H$ is a normal Abelian subgroup of $G$ such that 
$\omega|_{H \times H \times H}$ is cohomologically trivial and
$\mu \in (\Omega_{H, \omega})^G$ (defined in \eqref{Omega 1}). 

Theorem \ref{bij 1} showed that Lagrangian subcategories
of $\Rep(D^{\omega}(G))$ are in bijection with pairs $(H, \, B)$,
where $H$ is a normal Abelian subgroup of $G$ such that 
$\omega|_{H \times H \times H}$ is cohomologically trivial and
$B \in (\Lambda_{\omega}^2H)^G$ (the latter was defined in \eqref{alt omega}).

In this Subsection we will first show that the set of equivalence classes of indecomposable 
module categories over $\Vec_G^{\omega}$
such that the dual is pointed is in bijection with the set of Lagrangian subcategories
of $\Rep(D^{\omega}(G))$. 
Let $H$ be a normal Abelian subgroup of $G$ such that 
$\omega|_{H \times H \times H}$ is cohomologically trivial.
We will establish the aforementioned bijection by showing that there is a bijection 
between $\Omega_{H, \omega}$ (defined in \eqref{Omega}) 
and $\Lambda_{\omega}^2H$ that restricts to a bijection
between $(\Omega_{H, \omega})^G$ and $(\Lambda_{\omega}^2H)^G$.

Let $\mu \in C^2(H, \, k^\times)$ be a $2$-cochain satisfying
$\delta^2\mu = \omega|_{H \times H \times H}$. Define $alt^\prime(\mu)$ by
\[
alt^\prime(\mu)(h_1, \, h_2) := \frac{\mu(h_2, \, h_1)}
{\mu(h_1, \, h_2)}, h_1, h_2 \in H.
\]

\begin{lemma}
The map $alt^\prime(\mu): H \times H \to k^\times$ defined above
is an element of $\Lambda_{\omega}^2H$.
\end{lemma}
\begin{proof}
Clearly $alt^\prime(\mu)(h_1, \, h_2) = alt^\prime(\mu)(h_2, \, h_1)^{-1}$
and $alt^\prime(\mu)(h, \, h) = 1$, for all $h, h_1, h_2 \in H$.
We have
\begin{equation*}
\begin{split}
\frac{alt^\prime(\mu)(h, \, h_1) \,\,\, alt^\prime(\mu)(h, \, h_2)}
{alt^\prime(\mu)(h, \, h_1h_2)}
&= \frac{\mu(h_1, \, h)}{\mu(h, \, h_1)} 
\times \frac{\mu(h_2, \, h)}{\mu(h, \, h_2)}
\times \frac{\mu(h, \, h_1h_2)}{\mu(h_1h_2, \, h)}\\
&= \frac{\mu(h_1, \, h) \mu(h_2, \, h) \mu(hh_1, \, h_2)}
{\mu(h, \, h_2) \mu(h_1h_2, \, h) \mu(h_1, \, h_2)}
\times \omega(h, \, h_1, \, h_2)\\
&= \frac{\mu(h_1, \, h) \mu(hh_1, \, h_2)}
{\mu(h, \, h_2) \mu(h_1, \, hh_2)}
\times \omega(h, \, h_1, \, h_2)
\times
\omega(h_1, \, h_2, \, h)\\
&= \frac{\omega(h, \, h_1, \, h_2)\omega(h_1, \, h_2, \, h)}
{\omega(h_1, \, h, \, h_2)}\\
&= \beta_h(h_1, \, h_2),
\end{split}
\end{equation*}
\noindent for all $h, h_1, h_2 \in H$. In the second, third, and fourth
equalities above we used \eqref{mu omega} with
$(h_1, \, h_2, \, h_3) = (h, \, h_1, \, h_2),
(h_1, \, h_2, \, h), 
(h_1, \, h, \, h_2)$, respectively. 
\end{proof}

The map $alt^\prime$ induces a map between 
$\Omega_{H, \omega}$ and $\Lambda_{\omega}^2H$. By abuse of
notation we denote this map also by $alt^\prime$:
\begin{equation}
\label{alt prime}
alt^\prime:
\Omega_{H, \omega} \to \Lambda_{\omega}^2H :  
{\mu} \mapsto alt^\prime(\mu).
\end{equation}

\begin{lemma}
\label{bijec}
The map $alt^\prime$ defined above is a bijection.  
\end{lemma}
\begin{proof}
First note that $alt^\prime$ is well-defined.
Fix $\mu_0 \in C^2(H, \, k^\times)$ satisfying
$\delta^2\mu_0 = \omega|_{H \times H \times H}$.
Let $B_0 := alt^\prime(\mu_0)$.
Define bijections 
$f_1: \Lambda_{\omega}^2H \xrightarrow{\sim} \Lambda^2H : B \mapsto \frac{B}{B_0}$
and $f_2: \Omega_{H, \omega} \xrightarrow{\sim}
H^2(H, \, k^\times): {\mu} \mapsto {\lp\frac{\mu}{\mu_0}\rp}$.
Note that the cardinality of the two sets 
$\Omega_{H, \omega}$ and $\Lambda_{\omega}^2H$ are equal.
Injectivity, and hence bijectivity, of $alt^\prime$ follows from
the equality $f_1\circ alt^\prime = alt \circ f_2$.

%
%

\end{proof}

\begin{lemma}
\label{abc}
The following relation holds:
$$
\frac{\Upsilon_x(h_2, \, h_1)}{\Upsilon_x(h_1, \, h_2)} 
= \frac{\beta_{xh_1x^{-1}}(x, \, h_2)\beta_{xh_1x^{-1}}(xh_2, x^{-1})}
{\beta_{xh_1x^{-1}}(x, \, x^{-1})}, \mbox{ for all } x \in G, h_1, h_2 \in H.
$$
\end{lemma}
\begin{proof}
We have
\begin{equation*}
\begin{split}
&\frac{\Upsilon_x(h_2, \, h_1)}{\Upsilon_x(h_1, \, h_2)}
\times \frac {\beta_{xh_1x^{-1}}(x, \, x^{-1})}
{\beta_{xh_1x^{-1}}(x, \, h_2)\beta_{xh_1x^{-1}}(xh_2, x^{-1})} \\
&= \frac{\omega(xh_2x^{-1}, \, xh_1x^{-1},\, x)}
{\omega(xh_2x^{-1}, \, x, \, h_1) \omega(xh_1x^{-1}, \, xh_2x^{-1}, \, x)}
\times
\frac{\omega(xh_1x^{-1}, \, x, \, x^{-1}) \omega(x, \, x^{-1}, \, xh_1x^{-1})}
{\omega(x, \, h_1, \, x^{-1})}\\
& \hspace{3in} \times
\frac{\omega(xh_2, \, h_1, \, x^{-1})}
{\omega(xh_1x^{-1}, \, xh_2, \, x^{-1}) \omega(xh_2, \, x^{-1}, \, xh_1x^{-1})}\\
&= \frac {\omega(xh_2x^{-1}, \, xh_1x^{-1}, \, x) 
\omega(xh_1x^{-1}, \, x, \, x^{-1})\omega(x, \, x^{-1}, \, xh_1x^{-1})
\omega(xh_2x^{-1}, \, xh_1, \, x^{-1})} 
{\omega(xh_1x^{-1}, \, xh_2x^{-1}, \, x) \omega(xh_1x^{-1}, \, xh_2, \, x^{-1})
\omega(xh_2, \, x^{-1}, \, xh_1x^{-1}) \omega(xh_2x^{-1}, \, x, \, h_1x^{-1})}\\
&= \frac{\omega(x, \, x^{-1}, \, xh_1x^{-1}) 
\omega(xh_1h_2x^{-1}, \, x, \, x^{-1})}
{\omega(xh_1x^{-1}, \, xh_2x^{-1}, \, x)
\omega(xh_1x^{-1}, \, xh_2, \, x^{-1}) \omega(xh_2, \, x^{-1}, \, xh_1x^{-1})
\omega(xh_2x^{-1}, \, x, \, h_1x^{-1})}\\
& = \frac{\omega(x, \, x^{-1}, \, xh_1x^{-1}) \omega(xh_2x^{-1}, \, x, \, x^{-1})}
{\omega(xh_2, \, x^{-1}, \, xh_1x^{-1}) \omega(xh_2x^{-1}, \, x, \, h_1x^{-1})}\\
&= 1,
\end{split}
\end{equation*}
\noindent for all $x \in G, h_1, h_2 \in H$. In the first equality above we used the
definition of $\Upsilon$ and $\beta$ and canceled some factors. 
In the second, third, fourth, 
and fifth equalities we used \eqref{3-cocycle} with $(g_1, \, g_2, \, g_3, \, g_4)
= (xh_2x^{-1}, \, x, \, h_1, \, x^{-1}), \, (xh_2x^{-1}, \, xh_1x^{-1}, \, x, \, x^{-1}), \,
(xh_1x^{-1}, \, xh_2x^{-1}, \, x, \, x^{-1}),$ and 
$(xh_2x^{-1}, \, x, \, x^{-1}, \, xh_1x^{-1})$, respectively.
\end{proof}

\begin{lemma}
The map $alt^\prime$ defined in \eqref{alt prime} restricts to a 
bijection between $(\Omega_{H, \omega})^G$ and $(\Lambda_{\omega}^2H)^G$.
\end{lemma}
\begin{proof}
Let us first show that 
$alt^\prime((\Omega_{H, \omega})^G) \subseteq (\Lambda_{\omega}^2H)^G$. 
Pick any $\mu \in (\Omega_{H, \omega})^G$.\\
So $alt \lp \frac{\mu^x}{\mu} \times \Upsilon_x|_{H \times H} \rp = 1$,
for all $x \in G$. We have
\begin{equation*}
\begin{split}
alt^\prime(\mu)(x^{-1}ax, \, h) \times 
alt^\prime(\mu)(a, \, xhx^{-1})^{-1}
&= \frac{\mu(h, \, x^{-1}ax)}{\mu(x^{-1}ax, \, h)}
\times \frac{\mu(a, \, xhx^{-1})}{\mu(xhx^{-1}, \, a)}\\
&= \frac{\mu^x(x^{-1}ax, \, h)}{\mu(x^{-1}ax, \, h)}
\times \frac{\mu(h, \, x^{-1}ax)}{\mu^x(h, \, x^{-1}ax)}\\
&= alt \lp \frac{\mu^x}{\mu} \times \Upsilon_x|_{H \times H} \rp(h, \, x^{-1}ax)
\times \frac{\Upsilon_x(h, \, x^{-1}ax)}{\Upsilon_x(x^{-1}ax, \, h)}\\
&= \frac{\Upsilon_x(h, \, x^{-1}ax)}{\Upsilon_x(x^{-1}ax, \, h)}\\
&= \frac{\beta_a(x, \, h) \beta_a(xh, \, x^{-1})}{\beta_a(x, \, x^{-1})},
\end{split}
\end{equation*}
\noindent for all $x \in G, a \in H \cap R, h \in H$. In the fourth equality
above we used the fact that $alt \lp \frac{\mu^x}{\mu} \times \Upsilon_x|_{H \times H} \rp = 1$
and in the fifth equality we used Lemma \ref{abc}. So
$alt^\prime((\Omega_{H, \omega})^G)
\subseteq (\Lambda_{\omega}^2H)^G$, as desired. 

Now let us show that $(\Lambda_{\omega}^2H)^G 
\subseteq alt^\prime((\Omega_{H, \omega})^G)$.
Pick any $\mu \in \Omega_{H, \omega}$ and suppose that 
$alt^\prime(\mu) \in (\Lambda_{\omega}^2H)^G$. Suffices to show that 
$alt \lp \frac{\mu^x}{\mu} \times \Upsilon_x|_{H \times H} \rp = 1$,
for all $x \in G$. Let $B := alt^\prime(\mu)$. We have
\begin{equation*}
\begin{split}
&alt \lp \frac{\mu^x}{\mu} \times \Upsilon_x|_{H \times H} \rp(h_1, \, h_2)
\times \frac{\Upsilon_x(h_1, \, h_2)}{\Upsilon_x(h_2, \, h_1)}\\
& = B(xh_1x^{-1}, \, xh_2x^{-1}) B(h_1, \, h_2)^{-1}\\
&= B((yx^{-1})^{-1}a(yx^{-1}), \, xh_2x^{-1}) B(y^{-1}ay, \, h_2)^{-1} 
\qquad \qquad (\mbox{where $h_1 = y^{-1}ay$})\\
&= \frac{\beta_a(yx^{-1}, \, xh_2x^{-1}) \beta_a(yh_2x^{-1}, \,xy^{-1})}
{\beta_a(yx^{-1}, \, xy^{-1})}
\times
\frac{\beta_a(y, \, y^{-1})}{\beta_a(y, \, h_2) \beta(yh_2, \, y^{-1})}\\
&= \frac{\beta_a(yx^{-1}, \, xh_2) \beta_a(yh_2, \, x^{-1})
\beta_a(yh_2x^{-1}, \, xy^{-1})}
{\beta_{xh_1x^{-1}}(xh_2, \, x^{-1}) \beta_a(y, \, y^{-1})
\beta_a(yx^{-1}, \, xy^{-1}) \beta_a(y, \, h_2) \beta_a(yh_2, \, y^{-1})}\\
&= \frac{\beta_a(yx^{-1}, \, x) \beta_a(yh_2, \, x^{-1})
\beta_a(yh_2x^{-1}, \, xy^{-1}) \beta_a(y, \, y^{-1})} 
{\beta_{xh_1x^{-1}}(x, \, h_2) \beta_{xh_1x^{-1}}(xh_2, \, x^{-1})
\beta_a(yx^{-1}, \, xy^{-1}) \beta_a(yh_2, \, y^{-1})}\\
&= \frac{\beta_{xh_1x^{-1}}(x, \, x^{-1}) \beta_a(yh_2, \, x^{-1})
\beta_a(yh_2x^{-1}, \, xy^{-1}) \beta_a(y, \, y^{-1})}
{\beta_{xh_1x^{-1}}(x, \, h_2) \beta_{xh_1x^{-1}}(xh_2, \, x^{-1})
\beta_a(y, \, x^{-1}) \beta_a(yx^{-1}, \, xy^{-1})
\beta_a(yh_2, \, y^{-1})}\\
& = \frac{\Upsilon_x(h_1, \, h_2)}{\Upsilon_x(h_2, \, h_1)}
\times
\frac{\beta_a(y, \, y^{-1}) \beta_{h_1}(x^{-1}, \, xy^{-1})}
{\beta_a(y, \, x^{-1}) \beta_a(yx^{-1}, \, xy^{-1})}\\
& = \frac{\Upsilon_x(h_1, \, h_2)}{\Upsilon_x(h_2, \, h_1)},
\end{split}
\end{equation*}
\noindent for all $x \in G, h_1, h_2 \in H$. In the fourth through
eight equalities above we used \eqref{beta relation} with $(x, \, y, \, z)
= (yx^{-1}, \, xh_2, \, x^{-1}), (yx^{-1}, \, x, \, h_2), 
(yx^{-1}, \, x, \, x^{-1}),
(yh_2, \, x^{-1}, \, xy^{-1}),$ and $(y, \, x^{-1}, \, xy^{-1})$, respectively.
It follows that $(\Lambda_{\omega}^2H)^G 
\subseteq alt^\prime((\Omega_{H, \omega})^G)$ and
the Lemma is proved.
\end{proof}

Recall that $\mbox{E}(\C)$ denotes the set of (equivalence classes of)
braided tensor equivalences between a modular category $\C$ and
the centers of pointed fusion categories.

\begin{theorem}
\label{bijection 1}
The assignment $\M \mapsto \iota_\M$ $($defined in \eqref{iota}$)$ 
restricts to a bijection between equivalence classes of indecomposable 
$\Vec_G^\omega$-module categories  $\M$ with respect to which the 
dual fusion category $(\Vec_G^\omega)^*_\M$ is pointed and 
$E(\Rep(D^\omega(G)))$.
\end{theorem}
\begin{proof}
The proof is completely similar to the one given in Theorem \ref{bijection}.
\end{proof}

\begin{remark}
The equivalence type of the symmetric category $\L_{(H, B)}$ of
$\Rep(D^\omega(G))$ can be explicitly described in a way
similar to Remark~\ref{G'}, cf.\ \cite[Theorem 4.5]{N}.
\end{remark}

\begin{theorem}
\label{main 1}
Let $\C_1, \C_2$  be group-theoretical fusion categories. 
Then $\C_1,\C_2$ are weakly Morita equivalent if and only
if their centers $\Z(\C_1)$ and $\Z(\C_2)$ are equivalent as braided
fusion categories. 
\end{theorem}
\begin{proof}
That the  ``if" part is true for all fusion categories 
was first observed by M.~M\"uger in \cite[Remark 3.18]{Mu}.
This follows from the definition of weak Morita equivalence
and a theorem of P.~Schauenburg \cite{S}. See also \cite{Nat, O1, EO}.

For the ``only if" part, let $(G_1, \, \omega_1), (G_2, \, \omega_2)$ 
be two pairs of groups and $3$-cocycles such
that $\C_1$ is weakly Morita equivalent to $\Vec_{G_1}^{\omega_1}$
and $\C_2$ is weakly Morita equivalent to $\Vec_{G_2}^{\omega_2}$.

If $\Z(\C_1) \cong \Z(\C_2)$ (as braided fusion categories)
then  $\Z(\Vec_{G_1}^{\omega_1})\cong \Z(\Vec_{G_2}^{\omega_2})$
(as braided fusion categories) and therefore, $\Vec_{G_1}^{\omega_1}$
and $\Vec_{G_2}^{\omega_2}$ are weakly Morita equivalent by
Theorem~\ref{bijection 1} and hence, $\C_1$ and $\C_2$ are weakly
Morita equivalent.
\end{proof}

\begin{corollary}
\label{main 2}
Let $G, G'$ be finite groups, $\omega \in Z^3(G,\, k^\times)$, 
and $\omega' \in Z^3(G',\, k^\times)$. Then the representation
categories of twisted doubles $D^\omega(G)$ and 
$D^{\omega'}(G')$ are equivalent as braided tensor categories 
if and only if $G$ contains a normal Abelian subgroup $H$
such the following conditions are satisfied:
\begin{enumerate}
\item[(1)]  $\omega|_{H \times H \times H}$ is cohomologically trivial,
\item[(2)] there is a $G$-invariant $($see \eqref{Omega 1}$)$ 
$2$-cochain $\mu \in C^2(H, \, k^\times)$ such that
that $\delta^2 \mu = \omega|_{H \times H \times H}$, and
\item[(3)] there is
an isomorphism $a: G' \xrightarrow{\sim} \widehat{H} \rtimes_{\nu} 
(H \backslash G)$ such that $\varpi \circ (a \times a \times a)$ and 
$\omega'$ are  cohomologically equivalent. 
\end{enumerate}
Here $\nu$ is a certain $2$-cocycle in 
$Z^2(H \backslash G, \, \widehat{H})$ coming from the $G$-invariance
of $\mu$ and $\varpi$ is a certain $3$-cocycle on $\widehat{H} \rtimes_{\nu} 
(H \backslash G)$ depending on $\nu$ and on the exact sequence
$1 \to H \to G \to H \backslash G \to 1$ $($see \cite[Theorem 5.8]{N}
for precise definitions$)$.
\end{corollary}
\end{subsection}

\end{section}

\begin{section}{Examples}
\label{exs}

\subsection{Lagrangian subcategories of the Drinfeld doubles of finite symmetry groups}

\begin{example} 
\label{dihedral}
Consider the group of symmetries of a regular $n$-gon, i.e.,
the dihedral group 
$D_{2n}= \la r,\, s \mid r^n = s^2= 1,\, rs = sr^{-1} \ra ,\, n \geq 2$.
Let us describe the Lagrangian subcategories of  $\Rep(D(D_{2n})) \cong
\Z(\Vec_{D_{2n}})$.

Let $n= 2$. Then $D_4 \cong \mathbb{Z}/2\mathbb{Z} \times \mathbb{Z}/2\mathbb{Z}$.
There are six different  Lagrangian subcategories of  
$\Rep(D(\mathbb{Z}/2\mathbb{Z} \times \mathbb{Z}/2\mathbb{Z}))$,
all of them equivalent to 
$\Rep(\mathbb{Z}/2\mathbb{Z} \times \mathbb{Z}/2\mathbb{Z})$.

Let $n= 4$. The group $D_8$ has five normal Abelian subgroups which
give rise to seven Lagrangian subcategories of  $\Z(\Vec_{D_8})$.
With an exception of the Lagrangian subcategory corresponding to the center
$\la r^2\ra $ of $D_8$, all Lagrangian subcategories are equivalent to $\Rep(D_8)$.
The one corresponding to $\la r^2\ra $ is equivalent to 
$\Rep(\mathbb{Z}/2\mathbb{Z} \times \mathbb{Z}/2\mathbb{Z} \times \mathbb{Z}/2\mathbb{Z})$.
Applying \cite{DGNO} we conclude that 
for some $3$-cocycle $\omega$ on  $\mathbb{Z}/2\mathbb{Z} \times \mathbb{Z}/2\mathbb{Z} 
\times \mathbb{Z}/2\mathbb{Z}$ there is a braided tensor
equivalence $\Rep(D(D_8)) \cong
\Rep(D^\omega(\mathbb{Z}/2\mathbb{Z} \times \mathbb{Z}/2\mathbb{Z} \times \mathbb{Z}/2\mathbb{Z}))$
(existence of such an equivalence is already known to experts, see
\cite{CGR} and \cite{GMN}).

If $n =3$ or $n \geq 5$ then normal Abelian subgroups of $D_{2n}$ are precisely 
rotation subgroups
$\la r^k\ra $, where $k$ is a divisor of $n$. Each of these subgroups is cyclic and
so has a trivial Schur multiplier. The  Lagrangian
subcategory of  $\Z(\Vec_{D_{2n}})$ 
corresponding to $\la r^k\ra,\, k|n$, is
equivalent to $\Rep(\mbox{Dih}(\mathbb{Z}/\frac{n}{k}\mathbb{Z} \times
\mathbb{Z}/{k}\mathbb{Z}))$, where 
for any Abelian group $A$ we denote
$\mbox{Dih}(A) = A \rtimes \mathbb{Z}/2\mathbb{Z}$
the generalized dihedral group (the action of $ \mathbb{Z}/2\mathbb{Z}$ 
on $A$ is by inverting elements).
Consequently, there is a $3$-cocycle $\omega $ on 
$\mbox{Dih}(\mathbb{Z}/\frac{n}{k}\mathbb{Z} \times
\mathbb{Z}/{k}\mathbb{Z})$  such that
$\Rep(D^\omega(\mbox{Dih}(\mathbb{Z}/\frac{n}{k}\mathbb{Z} \times
\mathbb{Z}/{k}\mathbb{Z}))) \cong
\Rep(D(D_{2n}))$ as braided tensor categories.
\end{example}
Next,  consider the symmetry groups of Platonic solids.
\begin{example}
\label{platonic}
The tetrahedron group $A_4$ has two normal Abelian subgroups: the trivial one and
another isomorphic to $\mathbb{Z}/2\mathbb{Z} \times \mathbb{Z}/2\mathbb{Z}$. 
The Schur multiplier of $\mathbb{Z}/2\mathbb{Z} \times \mathbb{Z}/2\mathbb{Z}$
is isomorphic to $\mathbb{Z}/2\mathbb{Z}$ with the trivial $A_4$-action.
So  $\Rep(D(A_4))$ has three Lagrangian subcategories.
All three are equivalent to $\Rep(A_4)$.

The double of the cube/octahedron group $S_4$  also has  three Lagrangian 
subcategories (corresponding to the same data as in the case of $A_4$).
One can check that corresponding three Lagrangian subcategories
of $\Rep(D(S_4))$ are all isomorphic to $\Rep(S_4)$.

Finally, the group of symmetries of dodecahedron/icosahedron is
a simple non-Abelian group $A_5$. It is clear that for 
any simple non-Abelian group $G$ the category $\Rep(D(G))$
contains a unique Lagrangian subcategory (corresponding to the trivial
subgroup of $G$).
\end{example}

Recall \cite{N}  that a group $G$ is called 
{\em categorically Morita rigid} if
$\Vec_G$ being weakly Morita equivalent to $\Vec_{G'}^\omega$
implies  that groups $G$  and $G'$ are isomorphic. 
One can check  that  in this case $\omega$  must be  cohomologically 
trivial. Our results imply that  $G$ is categorically Morita rigid if
and only if all Lagrangian subcategories of $\Rep(D(G))$ are equivalent
to $\Rep(G)$ as symmetric categories.

It follows from Examples~\ref{dihedral} and \ref{platonic}
that groups $A_4,\, S_4,\, A_5$, as well as groups $D_{2n}$,
where $n$ is a square-free integer, are categorically Morita rigid.
It is clear that any group $G$ without non-trivial normal Abelian 
subgroups is categorically Morita rigid.

Finally, let us consider symmetries of vector spaces.

\begin{example}
Let $n$ be a positive integer and let $\mathbb{F}_q$ be the finite 
field with $q$ elements. Let $G= SL(n, \, q)$ denote the {\em special linear} 
group of $n \times n$ matrices with entries from $\mathbb{F}_q$, i.e.,
matrices having determinant equal to $1$. 

Let $PSL(n, \, q) = SL(n, \, q)/Z(SL(n, \, q))$,
where $Z(SL(n, \, q))$ is the center of $SL(n, \, q)$. 
Let us assume that $(n,\, q) \not = (2, \, 2), (2, \, 3)$.
It is  known that in this case $PSL(n, \, q)$ is simple.

Let $d= (n, \, q-1)$ be the greatest common divisor of $n$ and $q-1$.
The group $Z(SL(n, \, q)) \cong \mathbb{Z}/d\mathbb{Z}$ 
is cyclic and any normal subgroup
of $SL(n, \, q)$ is contained in $Z(SL(n, \, q))$. Thus,
Lagrangian subcategories of 
$\Rep(D(SL(n, \, q)))$ correspond to divisors of $d$.
One can easily describe the equivalence types of these
subcategories. For instance, the Lagrangian
subcategory corresponding to 
$Z(SL(n, \, q)) \cong \mathbb{Z}/d\mathbb{Z}$
is equivalent to 
$\Rep( \mathbb{Z}/d\mathbb{Z} \times PSL(n, \, q))$.
Therefore, $\Rep(D(SL(n, \, q))$ is equivalent (as a braided tensor category)
to $\Rep(D^{\omega}(\mathbb{Z}/d\mathbb{Z} \times PSL(n, \, q)))$
for some $3$-cocycle $\omega$.
\end{example}

\subsection{Non-pointed categories of dimension $8$}

\begin{example}
\label{eight}
It is known \cite{TY} that there are exactly four non-pointed 
fusion categories of dimension $8$ with integral dimensions
of objects: $\Rep(D_8)$; $\Rep(Q_8)$, where $Q_8$
is the group of quaternions; $KP$, the representation category
of the Kac-Paljutkin Hopf algebra \cite{KP}; and $TY$, the category
of representations of a unique $8$-dimensional quasi-Hopf algebra 
which is not gauge equivalent to a Hopf algebra \cite{TY}
(equivalently, $TY$ is the unique non-pointed fusion category 
of dimension $8$ with integral dimensions of objects which does not have a
fiber functor).

Let us show that these four categories belong to four different
weak Morita equivalence classes and hence, in view of
Theorem~\ref{main 1}, their centers are not equivalent as
braided tensor categories.

All proper subgroups of $Q_8$ are normal Abelian and have
trivial Schur multipliers. Hence all of them produce
pointed duals. So $\Rep(Q_8)$ 
is the only non-pointed dual of $\Vec_{Q_8}$.

The only non-normal subgroups of $D_8=\la r,\, s \mid r^4 = s^2= 1,\, rs = sr^3 \ra$ are reflection
subgroups $\la s\ra ,\, \la sr\ra$, and their conjugates. 
Such subgroups appear as factors in exact factorizations of
$D_8$ (one can take $\la r\ra $ as another factor). 
The corresponding dual categories $(\Vec_{D_8})^*_{\M(\la s\ra ,\, 1)}$
and $(\Vec_{D_8})^*_{\M(\la sr\ra ,\, 1)}$ admit fiber functors
by \cite[Corollary 3.1]{O2}
and so are representations of semisimple Hopf algebras.
But Hopf algebras corresponding to exact factorizations
of $D_8$ are known to be either commutative or cocommutative.
We conclude that
$$
(\Vec_{D_8})^*_{\M(\la s\ra ,\, 1)} \cong (\Vec_{D_8})^*_{\M(\la sr\ra ,\, 1)}
\cong \Rep(D_8),
$$
and, hence, $\Rep(D_8)$  is the unique non-pointed dual of $\Vec_{D_8}$.

Thus, neither $\Rep(Q_8)$  nor $\Rep(D_8)$  is weakly
Morita equivalent to any other non-pointed category.

It remains to check that the same is true for $TY$. 
Let $\omega_0$ be a non-trivial $3$-cocycle on
$D_8/ \la r^2,\, s\ra  \cong  \mathbb{Z}/2\mathbb{Z}$
(corresponding to the non-zero element of 
$H^3(\mathbb{Z}/2\mathbb{Z},\, k^\times) \cong \mathbb{Z}/2\mathbb{Z}$).
Let $\pi : D_8 \to D_8/ \la r^2,\, s\ra  $
be the canonical projection. Define a $3$-cocycle
$\omega$ on $D_8$  by
$\omega = \omega_0 \circ(\pi\times \pi\times \pi)$.
Then $\omega \equiv 1$ on $\la r^2,\, s\ra $ and 
the restrictions of $\omega$ on each of the subgroups
$\la r\ra $ and  $\la sr\ra $ are  cohomologically non-trivial.
This means that the complete list of equivalence classes of 
indecomposable module  categories over $\Vec_{D_8}^\omega$ 
consists of  $\M(\{ 1\},\, 1)$, $\M(\la r^2\ra ,\, 1)$, $\M(\la s\ra ,\, 1)$, 
and $\M(\la r^2,\, s\ra ,\, \mu)$, where $\mu\in
H^2(\mathbb{Z}/2\mathbb{Z}\times \mathbb{Z}/2\mathbb{Z},\,k^\times)$.
Therefore, the only non-pointed dual category of  $\Vec_{D_8}^\omega$ 
corresponds to  $\M(\la s\ra ,\, 1)$, see Example~\ref{pointed mod cats}.
It follows from the classification of fiber functors
on group-theoretical categories obtained 
in \cite[Corollary 3.1]{O2} that the category
$(\Vec_{D_8}^\omega)^*_{\M(\la s\ra ,\, 1)}$ does not have 
a fiber functor and hence 
$(\Vec_{D_8}^\omega)^*_{\M(\la s\ra ,\, 1)} \cong TY$.
Since all other duals of $\Vec_{D_8}^\omega$ are pointed,
it follows that $TY$ is not weakly Morita equivalent
to any other non-pointed fusion category.

Hence,  $\Rep(D_8)$, $\Rep(Q_8)$, $KP$, and $TY$
are pairwise weakly Morita non-equivalent fusion categories. 
Our claim about their centers follows from Theorem~\ref{main 1}.
\end{example}

Let us note that there is another $3$-cocycle $\eta$ on $D_8$
such that $(\Vec_{D_8}^\eta)^*_{\M(\la s\ra ,\, 1)} \cong KP$.
Up to a conjugation, such $\eta$ must have a trivial restriction 
on  $\la r^2,\, s\ra $ and $\la sr\ra $ (this can be seen from the Kac exact
sequence \cite{Kac}). 

\begin{remark}
The braided tensor equivalence classes of twisted doubles
of groups of order $8$ were studied in detail in \cite{GMN}
using the higher Frobenius-Schur indicators. In particular,
it was shown that there are precisely $20$ equivalence classes
of such { \em non-pointed} doubles. 
In view of results of the present paper this description can
be interpreted in terms of weak Morita equivalence classes 
of pointed categories of dimension $8$. 
\end{remark}

\end{section}
\bibliographystyle{ams-alpha}

\begin{thebibliography}{A} 

\bibitem[BK]{BK} B. Bakalov, A. Kirillov Jr., 
\textit{Lectures on tensor categories and modular functors},
University Lecture Series, vol. 21, American Mathematical Society,
Rhode Island, 2001.

\bibitem[CGR]{CGR} A. Coste, T. Gannon, and P. Ruelle,
\textit{Finite group modular data},
Nucl.Phys. \textbf{B581} (2000) 679-717.

\bibitem[De]{De}  P.~Deligne,
\textit{Cat\'egories tensorielles},
Moscow Math.\ J.\ \textbf{2} (2002), no. 2, 227-248. 

\bibitem[Dr]{D}  V.~Drinfeld,
\textit{On Poisson homogeneous spaces of Poisson-Lie groups}, 
Theoret. and Math. Phys. \textbf{95} (1993), no. 2, 524--525.

\bibitem[DGNO]{DGNO} V.~Drinfeld, S.~Gelaki, D.~Nikshych,
and V.~Ostrik,
\textit{Group-theoretical properties of nilpotent modular categories},
http://arxiv.org/math.QA/0704.0195.

\bibitem[DPR1]{DPR1} R. Dijkgraaf, V. Pasquier, and P. Roche,
\textit{Quasi-quantum groups related to orbifold models},
Nuclear Phys. B. Proc. Suppl. \textbf{18B} (1990), 60-72.

\bibitem[DPR2]{DPR2} R. Dijkgraaf, V. Pasquier, and P. Roche,
\textit{Quasi-Hopf algebras, group cohomology, and orbifold models},
Integrable systems and quantum groups (Pavia, 1990), 
World Sci.\ Publishing, River Edge, NJ, 75-98 (1992).

\bibitem[DVVV]{DVVV} R. Dijkgraaf, C.~Vafa, E.~Verlinde, and H.~Verlinde,
\textit{The operator algebra of orbifold models},
Commun.\ Math.\ Phys.\ \textbf{123} (1989), 485-526.


\bibitem[EG]{EG} P. Etingof and S. Gelaki,
\textit{Isocategorical Groups},
Int.\ Math.\ Res.\ Not., 2001, no. 2, 59-76.

\bibitem[EK]{EK} P.~Etingof, D.~Kazhdan,
\textit{Quantization of Poisson algebraic groups and Poisson homogeneous spaces},
Symetries quantiques (Les Houches, 1995), 935-946, North-Holland, Amsterdam, 1998. 

\bibitem[EM]{EM} S.~Eilenberg, S.~MacLane, 
\textit{Cohomology Theory of Abelian Groups and Homotopy theory I-IV},
Proc.\ Nat.\ Acad.\ Sci.\ USA,   \textbf{36} (1950), 443-447;
 \textbf{36}  (1950), 657-663; \textbf{37}  (1951), 307-310;
 \textbf{38}  (1952), 325-329.

\bibitem[ENO]{ENO} P. Etingof, D. Nikshych, and V. Ostrik,
\textit{On fusion categories},
Ann. of Math. \textbf{162} (2005), 581-642.

\bibitem[EO]{EO} P. Etingof and V. Ostrik,
\textit{Finite tensor categories},
Mosc. Math. J., 2003, no. 3, 627-654, 782-783.

\bibitem[FRS]{FRS} J.~Fuchs, I.~Runkel, and C.~Schweigert,
\textit{TFT construction of RCFT correlators III: Simple currents},
Nucl.\ Phys.\ B694 (2004) 277-353.

\bibitem[GMN]{GMN} C. Goff, G. Mason, and S. Ng,
\textit{On the Gauge Equivalence of Twisted Quantum Doubles of Elementary 
Abelian and Extra-Special 2-Groups}, 2006,
http://www.arxiv.org/abs/math.QA/0603191.

\bibitem[Kac]{Kac} G.I.~Kac, 
\textit{Extensions of groups to ring groups},
Math USSR Sb.\ \textbf{5} (1968), 451-474.

\bibitem[KacP]{KP} G.I.~Kac and V.G.~Palyutkin, 
\textit{Finite ring groups},
Trans.\ Moscow Math.\ Soc.\ \textbf{5} (1966), 251-294.

\bibitem[Ka]{KAR} G. Karpilovsky,
\textit{Projective representations of finite groups}, 
Monographs and textbooks in pure and applied mathematics,
Marcel Dekker, Inc, New York, 1985.

\bibitem[K]{K} C. Kassel,
\textit{Quantum groups},
Graduate texts in mathematics \textbf{155}, Springer-Verlag, New York, 1995.

\bibitem[MN]{MN}  G. Mason and S. Ng,
\textit{Group cohomology and gauge equivalence of some twisted quantum doubles},  
Trans. Amer. Math. Soc.,  \textbf{353}  (2001),  no. 9, 3465--3509 (electronic).

\bibitem[Mu1]{Mu} M. M\"uger,
\textit{From subfactors to categories and topology I. Frobenius algebras in 
and Morita equivalence of tensor categories},
J. Pure Appl. Algebra \textbf{180} (2003), 81-157.

\bibitem[Mu2]{M} M. M\"uger,
\textit{On the structure of modular categories},
Proc. London Math. Soc. (3) \textbf{87} (2003) 291-308.

\bibitem[Na]{N} D. Naidu,
\textit{Categorical Morita equivalence for group-theoretical categories},
Commun. Alg. (to appear). http://arxiv.org/math.QA/0605530.

\bibitem[N]{Nat} S. Natale,
\textit{On group theoretical Hopf algebras and exact factorization of
finite groups}, 
J.\ Algebra  \textbf{270}  (2003),  no. 1, 199--211. 

\bibitem[O1]{O1} V. Ostrik,
\textit{Module categories, weak Hopf algebras and modular invariants},
Transform. Groups \textbf{8} (2003), no.2, 177-206.

\bibitem[O2]{O2} V. Ostrik,
\textit{Module categories over the Drinfeld double of a finite group},
Int.\ Math.\ Res.\ Not., 2003, no. 27, 1507-1520.

\bibitem[Q]{Q} F. Quinn,
\textit{Group categories and their field theories},
Proceedings of the Kirbyfest (Berkeley, CA, 1998), 407-453, 
Geom. Topol. Monogr., 2, Geom. Topol. Publ., Conventry, 1999.

\bibitem [S]{S} P.~Schauenburg,
\textit{The monoidal center construction and bimodules},
J.\ Pure Appl.\ Alg.\ \textbf{158} (2001), 325-346.
 
\bibitem[TY]{TY} D.~Tambara, S.~Yamagami
\textit{Tensor categories with fusion rules of self-duality 
for finite abelian groups}, 
J.\ Algebra  \textbf{209}  (1998),  no. 2, 692--707.

\bibitem[T]{T} V.~Turaev, 
\textit{Quantum  invariants of knots and $3$-manifolds},
de Gruyter Stud.\ Math.\ Vol. 18, de Gruyter, Berlin 1994.

\end{thebibliography}

\end{document}